\documentclass{gtart_h}

\def\ifplaintex{\expandafter\ifx\csname documentclass\endcsname\relax}

\def\gtp{{\mathsurround=0pt\it $\cal G\mskip-2mu$eometry \&\ 
$\cal T\!\!$opology $\cal P\!$ublications}}  

\def\recd{{\small Received:\qua\receiveddate\ifx\reviseddate\relax
\else\qquad Revised:\qua\reviseddate\fi\par}} 

\def\lognumber#1{\def\thelognumber{#1}}
\def\volumenumber#1{\def\thevolumenumber{#1}}
\def\volumeyear#1{\def\thevolumeyear{#1}}
\def\papernumber#1{\def\thepapernumber{#1}}
\def\pagenumbers#1#2{\def\startpage{#1}\def\finishpage{#2}}
\def\published#1{\def\publishdate{#1}}

\def\received#1{\def\receiveddate{#1}}

\def\accepted#1{\def\accepteddate{#1}}

\def\asciiaddress#1{\def\theasciiaddress{#1}}
\def\asciiemail#1{\def\theasciiemail{#1}}

\long\def\asciiabstract#1{\long\def\theasciiabstract{#1}}

\let\\\par\let\thelognumber\relax\let\thevolumenumber\relax
\let\thepapernumber\relax\let\thevolumeyear\relax\let\startpage\relax
\let\finishpage\relax\let\publishdate\relax\let\receiveddate\relax
\let\reviseddate\relax\let\accepteddate\relax\let\theasciititle\relax
\let\theasciiauthors\relax\let\theasciiaddress\relax
\let\theasciiabstract\relax

\let\theasciiemail\relax

\ifplaintex
\font\logobig=cmssbx10 scaled 3836
\font\logomed=cmssbx10 scaled 2557
\else
\font\logobig=cmssbx10 scaled 4200
\font\logomed=cmssbx10 scaled 2800
\fi

\long\def\makeagttitle{   
\count0=\startpage
\agt\hfill      
\hbox to 45truept{\vbox to 0pt{\vglue -13truept{\logomed A\kern -.37em{\logobig 
T}\kern -.38em G}\vss}\hss}
\break
{\small Volume \thevolumenumber\ (\thevolumeyear)
\startpage--\finishpage\nl
Published: \publishdate}

\vglue .25truein

{\parskip=0pt\leftskip 0pt plus
1fil\def\\{\par\smallskip}{\Large\bf\thetitle}\par\medskip} \vglue
0.05truein

{\parskip=0pt\leftskip 0pt plus 1fil\def\\{\par}{\sc\theauthors}
\par\medskip}%
 
\vglue 0.03truein 

{\small\leftskip 25truept\rightskip 25truept{\bf Abstract}\stdspace\theabstract

{\bf AMS Classification}\stdspace\theprimaryclass
\ifx\thesecondaryclass\relax\else; \thesecondaryclass\fi\par
{\bf Keywords}\stdspace \thekeywords\par}\vglue 7truept

}   

\ifplaintex
\font\phead=cmsl9 scaled 950
\font\pnum=cmbx10 scaled 913
\font\pfoot=cmsl9 scaled 950
\headline{\vbox to 0pt{\vskip -4.5mm\line{\small\phead\ifnum
\count0=\startpage ISSN 1472-2739 (on-line) 1472-2747 (printed)
\hfill {\pnum\folio}\else\ifodd\count0\def\\{ }%
\ifx\theshorttitle\relax\thetitle\else\theshorttitle\fi\hfill{\pnum\folio}
\else\def\\{ and }{\pnum\folio}\hfill\ifx\theshortauthors\relax\theauthors
\else\theshortauthors\fi\fi\fi}\vss}}
\footline{\vbox to 0pt{\vglue 0mm\line{\small\pfoot\ifnum\count0=\startpage
\copyright\ \gtp\hfill\else
\agt, Volume \thevolumenumber\ (\thevolumeyear)\hfill\fi}\vss}}
\else
\font=cmsl9 scaled 1050
\font\lnum=cmbx10 
\font=cmsl9 scaled 1050
\makeatletter
\def\@oddhead{{\small\lhead\ifnum\count0=\startpage ISSN 1472-2739 
(on-line) 1472-2747 (printed)\hfill {\lnum\number\count0}\else\ifodd\count0
\def\\{ }\ifx\theshorttitle\relax \thetitle \else\theshorttitle\fi\hfill
{\lnum\number\count0}\else\def\\{ and }{\lnum\number\count0}
\hfill\ifx\theshortauthors\relax 
\theauthors\else\theshortauthors\fi\fi\fi}}\def\@evenhead{\@oddhead}
\def\@oddfoot{\small\lfoot\ifnum\count0=\startpage\copyright\ \gtp\hfill\else
\agt, Volume \thevolumenumber\ (\thevolumeyear)\hfill\fi}
\def\@evenfoot{\@oddfoot}
\makeatother
\fi
\let\maketitlepage\makeagttitle
\let\makeshorttitle\maketitlepage
\let\maketitle\maketitlepage

\newwrite\gtoutfile
\long\gdef\makeheadfile{  
{\def\\{, }\def\s{ }
\immediate\openout\gtoutfile head.xxx
\immediate\write\gtoutfile{Proxy-for: \ifx\theasciiauthors\relax
\theauthors\else\theasciiauthors\fi\s<\ifx\theasciiemail\relax\theemail\else\theasciiemail\fi>}
\immediate\write\gtoutfile{\noexpand\\}
\immediate\write\gtoutfile{Authors: \ifx\theasciiauthors\relax
\theauthors\else\theasciiauthors\fi}
{\def\\{ }\immediate\write\gtoutfile{Title: \ifx\theasciititle\relax
\thetitle\else\theasciititle\fi}}
\immediate\write\gtoutfile{Subj-class: GT or SG, GR etc}
\immediate\write\gtoutfile{MSC-class: \theprimaryclass\ifx\thesecondaryclass\relax\else, \thesecondaryclass\fi}
\immediate\write\gtoutfile{Journal-ref: Algebr. Geom. Topol. \thevolumenumber\s
(\thevolumeyear) \startpage-\finishpage}
\immediate\write\gtoutfile{Comments: Published by Algebraic and
Geometric Topology at}
\immediate\write\gtoutfile{\s\s\s  http://www.maths.warwick.ac.uk/agt/AGTVol\thevolumenumber/agt-\thevolumenumber-\thepapernumber.abs.html}
\immediate\write\gtoutfile{\noexpand\\}
\immediate\write\gtoutfile{}
\ifx\theasciiabstract\relax
\immediate\write\gtoutfile{\theabstract}\else
\immediate\write\gtoutfile{\theasciiabstract}\fi
\immediate\write\gtoutfile{}
\immediate\write\gtoutfile{\noexpand\\}
\immediate\write\gtoutfile{}
\immediate\closeout\gtoutfile}}  

\def\maketitlepage{\makeagttitle\makeheadfile}
\let\makeshorttitle\maketitlepage
\let\maketitle\maketitlepage

\lognumber{22}
\volumenumber{4}
\volumeyear{2004}
\papernumber{22}
\published{27 June 2004}
\pagenumbers{439}{472}
\received{10 April 2003}
\accepted{20 May 2004}

\usepackage{amssymb,amsmath,graphicx}

\newtheorem{thm}{Theorem}
\newtheorem{prop}[thm]{Proposition}
\newtheorem{lemma}[thm]{Lemma}
\newtheorem{cor}[thm]{Corollary}

\theoremstyle{definition}
\newtheorem{defn}[thm]{Definition}
\newtheorem{rem}[thm]{Remark}
\newtheorem{exam}[thm]{Example}

\def\Cal#1{{\cal#1}}
\def\<{\langle}\def\>{\rangle}\def\ov{\overline}
\def\wtil{\widetilde}
\def\Z{{\mathbb Z}}
\def\N{{\mathbb N}} 
\def\R{{\mathbb R}}
\def\E{{\mathbb E}}

\def\ep{\epsilon}               \def\varep{\varepsilon}
\def\Remark{\paragraph{Remark}}
\let\noproof\qed

\begin{document}

\title{Embeddings of graph braid and surface
groups\\in right-angled Artin groups and braid groups}
\shorttitle{Embeddings in right-angled Artin groups and braid groups}

\authors{John Crisp\\Bert Wiest}
\address{Institut de Math\'ematiques de Bourgogne (IMB), UMR 5584 du CNRS\\
Universit\'e de Bourgogne, 9 avenue Alain Savary, B.P. 47870\\ 
21078 Dijon cedex, France\\
\smallskip
IRMAR, UMR 6625 du CNRS, Campus de Beaulieu, Universit\'e de Rennes 1\\
35042 Rennes, France}
\asciiaddress{Institut de Mathematiques de Bourgogne (IMB), UMR 5584 du CNRS\\
Universite de Bourgogne, 9 avenue Alain Savary, B.P. 47870\\ 
21078 Dijon cedex, France\\
IRMAR, UMR 6625 du CNRS, Campus de Beaulieu, Universite de Rennes 1\\
35042 Rennes, France}
\gtemail{\mailto{jcrisp@u-bourgogne.fr}, \mailto{bertold.wiest@math.univ-rennes1.fr}}
\asciiemail{jcrisp@u-bourgogne.fr, bertold.wiest@math.univ-rennes1.fr}

\begin{abstract}
We prove by explicit construction that graph braid groups and most surface
groups can be embedded in a natural way in right-angled Artin groups, and we 
point out some consequences of these embedding results. We also show that every
right-angled Artin group can be embedded in a pure surface braid group.
On the other hand, by generalising to right-angled Artin groups a result
of Lyndon for free groups, we show that the Euler characteristic $-1$ surface
group (given by the relation $x^2y^2=z^2$) never embeds in a right-angled
Artin group.
\end{abstract}
\asciiabstract{%
We prove by explicit construction that graph braid groups and most surface
groups can be embedded in a natural way in right-angled Artin groups, and we 
point out some consequences of these embedding results. We also show that every
right-angled Artin group can be embedded in a pure surface braid group.
On the other hand, by generalising to right-angled Artin groups a result
of Lyndon for free groups, we show that the Euler characteristic -1 surface
group (given by the relation x^2y^2=z^2) never embeds in a right-angled
Artin group.}


\maketitle 


\section{Introduction and statement of results}


Right-angled Artin groups (which also go by the names of graph groups or
partially commutative groups) are by definition those finitely presented groups
whose defining relations consist simply of a certain number of commutation relations 
$xy=yx$ between elements $x$, $y$ of the generating set. It is customary
and convenient to specify a right-angled Artin group by drawing a graph
whose vertices correspond to the generators, and whose edges connect
``commuting'' vertices. These groups interpolate between the finite rank free
and free abelian groups and share a number of important
properties.  They are linear (so residually finite) \cite{Hum,HW},
torsion free, and even bi-orderable as well as residually nilpotent \cite{DT}.
Moreover, right-angled Artin groups act freely properly
co-compactly on certain naturally associated CAT(0) cubed complexes, which
implies that they are biautomatic (by~\cite{NR}, see also \cite{VWyk,HM})
and have quadratic isoperimetric inequality. Even better, these groups admit
a polynomial time solution to the conjugacy problem \cite{Serv,Wrath}.
Some but not all of the above properties are inherited by subgroups.
On the other hand, subgroups of right-angled Artin groups 
can exhibit some surprising properties, depending on the structure 
of the underlying defining graph. 
Bestvina and Brady \cite{BB} were able, for instance,
to distinguish for the first time between a variety
of  subtly different finiteness properties by looking at certain subgroups
of these groups. 
Graph groups and their subgroups are therefore of considerable
interest, not only for their  attractive properties, but also as a varied
source of examples.

In Section \ref{Sect2} we shall describe briefly the geometric 
properties of cubed complexes.
We shall restate Gromov's combinatorial characterisation of the 
CAT(0) property for these complexes 
and shall prove a result which allows one to identify locally 
isometric embeddings between such spaces.
This gives a method for constructing (quasi-convex) embeddings 
between groups which are fundamental groups of
locally CAT(0) cubed complexes.

There are two classes of groups that we shall embed in right-angled Artin groups,
using this technique:  graph braid groups 
(i.e.~fundamental groups of configuration spaces of $k$ disjoint points in a graph)
which are treated in Section \ref{Sect3}, and surface groups, treated in 
Section \ref{Sect4}.  All of these examples with the exception of exactly three
closed surface groups can be successfully embedded.

When discussing graph braid groups it is more convenient to consider
the larger class of
``reduced graph braid groups'' (defined in Section \ref{Sect3}). These are 
the fundamental
groups of certain cubed complexes which, unlike the usual configuration spaces,
depend on the simplicial structure of the underlying graph. However,
the braid group $B_n(G)$
over a graph $G$ is always isomorphic to the reduced graph braid group $RB_n(G')$
for a sufficiently fine subdivision $G'$ of $G$ (see \cite{Ab1,Ab2,AbGh}).
It was known already from
\cite{Ab1,AbGh} that reduced graph braid groups act freely co-compactly
on CAT(0) cubed complexes, and so are biautomatic, by the result of \cite{NR}.
In Section \ref{Sect3} we present a very natural
embedding of each graph braid group in a right-angled Artin group which is
induced by a locally isometric embedding of locally CAT(0) cubed
complexes. As corollaries we have that graph braid groups are linear,
bi-orderable, residually finite, residually nilpotent. 
In a future article we plan to give a proof that
these groups also admit a polynomial time
solution to the conjugacy problem.  

There are two particular examples of (reduced) graph braid groups that are
isomorphic to surface groups. These are the 2-string braid groups
$B_2(K_{3,3})\cong RB_2(K_{3,3})$ and $B_2(K_5)\cong RB_2(K_5)$,
where $K_{3,3}$ denotes the utilities
graph, and $K_5$ the complete graph on $5$ vertices (the smallest two
non-planar graphs). The associated cubed complexes are squarings of
non-orientable surfaces of Euler characteristic $-3$ and $-5$
respectively. See \cite{Ab1} for a discussion of these examples.
Our construction embeds $B_2(K_5)$ into the right-angled
Artin group whose defining graph is the 1-skeleton of the twofold quotient
of the dodecahedron by the antipodal map.  The group $B_2(K_{3,3})$
embeds in the right-angled Artin group with defining graph
$\Delta(K_{3,3})$ shown in Figure \ref{F:chi-3}.

Results of \cite{DSS} have shown that `most' right-angled Artin groups contain 
non-abelian closed surface subgroups. More precisely, if the defining graph
contains an induced (full) subgraph isomorphic to an $n$-cycle for $n\geq 5$, 
then the right-angled Artin group contains a non-abelian closed surface subgroup. 
 However all the examples from
\cite{DSS} are orientable surface groups of relatively high genus
(at least $5$), which raises the question of exactly which surface groups can 
be embedded in a right-angled Artin group. This question is answered in Section
\ref{Sect4} where we show that all surface groups 
(orientable, nonorientable, closed or with boundary) can
be embedded in some right-angled Artin group, with three exceptions:
the fundamental groups of the 
projective plane (which has 2-torsion), the Klein bottle (which has generalised
torsion, so is not bi-orderable), and the nonorientable surface $S_{-1}$ with
Euler characteristic $-1$. We remark that the dual question of classifying those 
right-angled Artin groups (or Artin groups) which contain non-abelian
(hyperbolic) closed surface subgroups remains open.
Some partial results, notably in the case of finite and 
affine type Artin groups, are given in \cite{GLR}.

The group $\pi_1(S_{-1})$, the last of the three exceptional surface groups
mentioned above, has presentation $\<x,y,z\mid x^2y^2=z^2\>$. In order to
prove that this group cannot embed in any right-angled Artin group
 we give a topological proof
that if elements $x,y,z$ of a right-angled Artin group satisfy the
relation $x^2y^2=z^2$ then $x$, $y$ and $z$ mutually commute (see Theorem
\ref{x2y2z2}). This generalizes a classic result of Lyndon on free groups
\cite{Lyn} to the case of right-angled Artin groups.  Lyndon's result was
shortly afterwards generalized to the following statement \cite{LS}: 
in a free group
$F$, the relation $x^my^n=z^k$ implies that $x,y$ and $z$ commute.  We do
not know whether this statement also holds true over an arbitrary
right-angled Artin group.

In Section \ref{Sect5} we use somewhat different techniques to those
above, to prove a result about  embeddings in \emph{braid groups}.
We show that every right-angled Artin group (of rank $n$) may be
embedded as a subgroup of the pure braid group
(on $n$ strings) over some compact surface with boundary. Often this
surface can be chosen to be planar, and so the right-angled Artin group
can be embedded in a classical pure braid group $B_l$ for some $l\in\N$.
This is discussed in Section \ref{Sect6}.
However, we are not able to answer the following question: does every
right-angled Artin group embed in a classical braid group? We also leave
open the question whether the surface group $\<x,y,z\mid x^2y^2=z^2\>$
embeds in a braid group $B_l$ for some $l\in\N$.

{\bf Acknowledgement}\qua We wish to thank Aaron Abrams for 
stimulating discussions in the early stages of this work.
We are also grateful to Mark Sapir and Ilya Kapovich
for helpful remarks on an earlier version of this manuscript,
and to Luis Paris for bringing to our attention the work of Duchamp and Thibon on 
orderability and related properties of right-angled Artin groups.
Part of the research was done while Bert Wiest
was a postdoc at the Pacific Institute for the Mathematical Sciences (PIMS).


\section{Background on CAT(0) cube complexes}\label{Sect2}


By a \emph{cubed complex} $X$ we shall mean a polyhedral complex in which
each cell is a finite dimensional unit Euclidean cube.  We refer to
\cite{BH}, pp.114--115, for a precise definition of polyhedral complex.  
Briefly, a cubed complex $X$ is an identification space defined by an
equivalence relation on the disjoint union of a collection of Euclidean
cubes, wherein the interior of each face of each cube is mapped
injectively to $X$, and if two points $x,x'$ lying in the interior of
faces $F,F'$ are identified then there is an isometry $h\co F\to F'$ such
that $y$ and $h(y)$ are identified for all $y\in F$. For example, the
$n$-torus obtained by identifying opposite faces of an $n$-cube is a cubed
complex.

Every finite dimensional cubed complex $X$ is equipped with its so-called
\emph{intrinsic metric} $d$, defined by setting $d(x,y)$ to be the infimum
of lengths of rectifiable paths from $x$ to $y$ in $X$ (where a
rectifiable path is a path obtained by concatenating the images of
finitely many straight line segments in cubes and the length of such a
path is the sum of the Euclidean lengths of its segments). This actually
defines a complete length metric on $X$ (see \cite{BH}).

By a {\sl CAT(0)} {\it space} we mean a geodesic metric space which
satisfies the CAT(0) triangle inequality as defined, for example, in
\cite{Gro,BH}.  This is a global nonpositive curvature condition which
states that distances between points on a geodesic triangle are no
greater than the distances between the corresponding points on a
reference triangle of the same sidelengths in the Eulidean plane.  By
a \emph{locally} {\sl CAT(0)} {\it space} we mean metric space in
which every point has an open ball neighbourhood which is a CAT(0)
space with the induced metric. The Cartan-Hadamard Theorem
\cite{Gro,BH} states that if $X$ is a complete locally CAT(0) space
then its universal cover $\wtil X$ is a complete CAT(0) space.

By a \emph{local isometry} between metric spaces we mean a map $f\co X\to
Y$ such that every point $p$ in $X$ has an open ball neighbourhood $N_p$
in $X$ such the restriction of $f$ to $N_p$ is an isometry with respect to
the induced metric on $N_p$. If $X$ is a geodesic metric space, $Y$ a
locally CAT(0) space, and $f\co X\to Y$ a local isometry, then $X$ is also
locally CAT(0) and $f$ is $\pi_1$-injective. In fact $f$ lifts to an
isometric embedding of $\wtil X$ in $\wtil Y$ (see \cite{BH}, p.201).

Any combinatorial map between cubed complexes which takes each cube to a
cube of the same dimension while respecting the face relation determines
uniquely a map between the metric cubed complexes which is a local
isometry on the interior of each cube. We call such a map a \emph{cubical
map}.  The following Theorem will be useful in proving the existence of
$\pi_1$-injective maps between cubed complexes. By a \emph{flag complex}
we mean a simplicial complex $K$ in which every set of $n+1$ vertices
which spans a complete graph in the $1$-skeleton of $K$ spans an
$n$-simplex in $K$.

Part (1) of the following Theorem is due to Gromov \cite{Gro}, and is
widely known. Part (2) seems to be the natural extension of Gromov's
idea, and is essentially a consequence of Lemmas 1.4--1.6 of
\cite{Ch}.
 
\begin{thm}\label{GromovThm}

{\rm(1)}\qua A finite dimensional cubed complex $Y$ is locally CAT(0) if and only
the link of every vertex is a flag complex (or equivalently, if and only
if, for every cube $C$ in $Y$, every triangle in the 1-skeleton of
$Lk(C,Y)$ bounds a 2-simplex).

{\rm(2)}\qua Let $X$ and $Y$ be finite dimensional cubed complexes and $\Phi\co
X\to Y$ a cubical map. Suppose that $Y$ is locally CAT(0). Then the map
$\Phi$ is a local isometry if and only if, for every vertex $x\in X$, the
simplicial map $Lk(x,X)\to Lk(\Phi(x),Y)$ induced by $\Phi$ is injective
with image a full subcomplex of $Lk(\Phi(x),Y)$.

Moreover, in this case, $X$ is locally CAT(0) and $\Phi$ is
$\pi_1$-injective and lifts to an isometric embedding of $\wtil X$ in
$\wtil Y$. 
\end{thm}

\Remark
If, in part (2) of the Theorem, it is already known that the links in $X$
are flag (equivalently that $X$ is locally CAT(0)) then it is enough to
check that the $1$-skeleton of $Lk(x,X)$ is mapped isomorphically onto a
full subgraph of the $1$-skeleton of $Lk(\Phi(x),Y)$, for all vertices
$x\in X$, in order to show that $\Phi$ is a local isometry.

\begin{proof}
We refer to \cite{Gro} (see also \cite{BH}) for the proof of (1). Note
that the link of any vertex in a cube complex is an \emph{all-right
spherical complex}, a piecewise spherical simplicial complex in which all
edge lengths and all dihedral angles measure $\frac{\pi}{2}$.  Note also
that, in an all-right spherical complex $L$, the open ball $B(v,\frac{\pi}{2})$
of radius $\frac{\pi}{2}$ around a vertex $v$ consists exactly of $v$ and the
interior of every simplex which has $v$ as a vertex. The sphere
$S(v,\frac{\pi}{2})$ of radius $\frac{\pi}{2}$ is simply the simplicial
link of $v$ in $L$. The following observation is a key element in Gromov's
proof of (1): any locally geodesic segment in $L$ which has endpoints in
$S(v,\frac{\pi}{2})$ and intersects $B(v,\frac{\pi}{2})$ has length at
least $\pi$. We will use this to prove part (2).

It is easily verified that the condition given in (2) is necessary for
$\Phi$ to be a local isometry. To prove sufficiency we use the fact
(an easy consequence of Lemmas 1.4 and 1.5 of \cite{Ch}) that a
cubical map $\Phi\co X\to Y$ is locally isometric if, for every vertex
$x\in X$, the induced map $\Phi_x\co Lk(x,X)\to Lk(\Phi(x),Y)$ is
$\pi$-distance preserving. A map $f\co L'\to L$ between piecewise
spherical complexes is said to be \emph{$\pi$-distance preserving}, if
$d_{L'}(p,q)\geq\pi$ implies $d_L(f(p),f(q))\geq\pi$, for all $p,q\in
L'$.

Suppose that $x$ is a vertex of $X$ and that $\Phi_x$ is injective
with image a full subcomplex $K$ of $L=Lk(\Phi(x),Y)$. We show that
$\Phi_x$ is $\pi$-distance preserving. The argument follows that of
Lemma 1.6 of \cite{Ch}. Given $p,q\in K$ such that $d_L(p,q)<\pi$
there exists a geodesic segment $\gamma\in L$ joining $p$ to
$q$. (Note that links in a CAT(0) polyhedral complex are
$\pi$-geodesic \cite{BH}). It is enough to show that $\gamma$ lies
wholly in $K$.  Suppose that $\gamma$ enters the interior of some
simplex $\sigma$ of $L$ and let $v$ denote a vertex of $\sigma$. Then
$\gamma$ enters $B(v,\frac{\pi}{2})$. Now, by the observation made
above, $\gamma$ must have an endpoint $p$ in $B(v,\frac{\pi}{2})$
(since otherwise it would contain a subsegment of length at least
$\pi$). Now any simplex which contains $p$ has $v$ as a vertex, and
since $p$ lies in the subcomplex $K$, we must have $v\in K$.  Since
this is true for every vertex of $\sigma$, and $K$ is a full
subcomplex, we deduce that $\sigma$ is a simplex of $K$.  But that is
to say that $\gamma$ never enters the interior of a simplex which does
not already lie in $K$, and so $\gamma$ lies wholly in $K$.
\end{proof}


\section{Graph braid groups and right-angled Artin groups}\label{Sect3}


Let $G$ be a finite graph, viewed as a simplicial complex. For each point
$x\in G$ denote by $c(x)$ the \emph{carrier} of $x$, that is the smallest
dimensional simplex which contains $x$. Thus $c(x)$ denotes either a
vertex or a (closed) edge of $G$. For $n\in\N$, we define $\wtil X_n(G)$
to be the subspace of $G^n$ consisting of those $n$-tuples
$(x_1,x_2,\ldots,x_n)$ for which the $x_i$ have mutually disjoint carriers:
$c(x_i)\cap c(x_j)=\emptyset$ for $i\neq j$. This is, in fact, a subspace
of the ordered $n$-point configuration space of $G$. The symmetric group
$S_n$ acts freely and properly discontinuously on $\wtil X_n(G)$ by
permutations of the coordinates, and the quotient space may be viewed as a
subspace of the unordered $n$-point configuration space of $G$ and is
written
\[
X_n(G)=\wtil X_n(G)/S_n =
\{ \{x_1,x_2,\ldots,x_n\} : c(x_i)\cap c(x_j)=\emptyset \text{ for } i\neq j\}\,.
\]

Following Abrams \cite{Ab1,Ab2} we refer to $X_n(G)$ (respectively 
$\wtil X_n(G)$) as the \emph{reduced} unordered (respectively ordered)
$n$-point configuration space of $G$. We concern ourselves here with the
\emph{reduced braid group} of $G$, namely \[RB_n(G)=\pi_1(X_n)\,.\] (This
group contains the \emph{reduced pure braid group} $RPB_n(G)=\pi_1(\wtil
X_n(G))$ as an index $n!$ subgroup).

Note that these definitions depend upon the simplicial structure on $G$,
while the usual configuration spaces of $G$ do not. It was proved in
\cite{Ab2,AbGh} that, after subdividing the edges of the graph $G$
sufficiently often, creating some valence-two vertices, the reduced
configuration space $X_n(G)$ becomes a deformation retract of the usual
configuration space of $G$, and hence that the reduced braid and reduced
pure braid groups coincide with the usual braid and pure braid groups,
$B_n(G)$ and $PB_n(G)$, respectively.  Thus the reader should keep in mind
that our results about $RB_n(G)$ hold in particular for the full graph
braid groups $B_n(G)$.

The reduced configuration space $X_n(G)$ has the structure of a cubed
complex, in a very natural way. The $k$-dimensional cubes are in 1-1
correspondence with collections $\{ c_1,c_2,\ldots,c_n\}$ of mutually disjoint
simplices in $G$, $k$ of which are edges and $n-k$ of which are vertices.
The cube associated to $\{ c_1,c_2,\ldots,c_n\}$ is simply the set of points
\[
\Cal C\{c_1,c_2,\ldots,c_n\}=\{ \{x_1,x_2,\ldots,x_n\} : x_i\in c_i 
\text{ for all } i\} \cong c_1\times c_2\times\cdots\times c_n\,.
\]
In fact, it is easily checked that $X_n(G)$ is a \emph{cubical} complex
(in the language of \cite{BH}) -- each cube is embedded and the intersection
of any two cubes is either empty or a single face of some dimension. 
It is also an almost immediate consequence of the
above description of the cubical structure on $X_n(G)$ that all vertex
links are flag complexes, and hence that $X_n(G)$ is a locally CAT(0)
cubical complex as first proved by Abrams \cite{Ab2,Ab1}. In this paper we
construct a locally isometric cubical map from $X_n(G)$ to a known CAT(0)
cubed complex whose fundamental group is a right-angled Artin group.  In
this way, we prove that every reduced graph braid group embeds in a
right-angled Artin group.

\paragraph{Definition} (Right-angled Artin group/graph group)\qua
Let $\Delta$ denote a finite simple graph (that is a finite graph with no
loops or multiple edges).
Let $\Cal V$ and $\Cal E$ denote the vertex and edge sets of $\Delta$
respectively. Associated to each such graph $\Delta$ is the \emph{right
angled Artin group}, or \emph{graph group}
\[ A_\Delta =\<\ \Cal V\ |\  uv=vu\ \text{ if } (u,v)\in\Cal E\ \>\,. \] 

We shall now describe a locally CAT(0) cubed complex which is an
Eilenberg-MacLane complex for $A_\Delta$. This construction is well-known,
and plays a crucial role for instance in the paper of Bestvina and Brady
\cite{BB} on finiteness properties of groups.

The standard $n$-torus $T^n=\E^n/\Z^n$ may be described as the unit
Euclidean cube $[0,1]^n$ with opposite faces identified, and so is a cubed
complex with one vertex, $n$ edges and precisely one $k$-cube spanned by
each set of $k$ distinct edges.  To each finite set $S$ of size $n$ we may
associate such a cubed $n$-torus $T(S)$ by labelling the edges by the
elements of $S$. We shall also orient all edges.

Given the graph $\Delta$ we now define the cubed complex $\Cal T_\Delta$
which is obtained from $T(\Cal V)$ by deleting the interior of every cube
which correspond to a set of vertices in $\Delta$ which does not span a
complete graph. Alternatively,
\[
\Cal T_\Delta =
\left( \coprod \{T(U) : U\subset\Cal V \text{ spans a complete graph in }
\Delta\} \right )/\sim
\]
where the $\sim$ denotes identification of $T(Q)$ as a subcomplex of
$T(R)$ in the obvious way whenever $Q\subset R$. Note that $\Cal
T(\Delta)$ has a unique vertex, the ``0-cube" $T(\emptyset)$. The link of
this vertex, which we shall denote $L_\Delta$, has two vertices $v^+$ and
$v^-$ for each $v\in\Cal V$ (corresponding to where the edge $v$ leaves
and enters the vertex). A set
$\{v_0^{\varep_0},v_1^{\varep_1},\ldots,v_k^{\varep_k}\}$ of $k$ distinct
vertices (with $\varep_i\in\{\pm\}$) spans a $k$-simplex in $L_\Delta$ if
and only if $v_0,\ldots,v_k$ span a complete graph in $\Delta$.

It is now obvious that $L_\Delta$ is a flag complex, and so by Theorem
\ref{GromovThm}(1) that $\Cal T_\Delta$ is a locally CAT(0) cubed complex.  
Moreover, it is clear that $\pi_1(\Cal T_\Delta)\simeq A_\Delta$ where the
standard generators of $A_\Delta$ are simply represented by the
corresponding labelled oriented edges of $\Cal T_\Delta$.

\medskip

We return now to our finite graph $G$, and define a new graph $\Delta(G)$
whose vertex set is just the set 
$\{ \overline x : x \text{ an edge of }G\}$ 
and where two vertices are joined by an edge in $\Delta(G)$ whenever the
corresponding edges in $G$ are disjoint (as closed sets). This definition
immediately allows the definition of a cubical map $\Phi \co X_n(G) \to
\Cal T_{\Delta(G)}$ as follows. Choose an orientation for each edge in
$G$, which naturally induces an orientation on each edge of $X_n(G)$.
Recall that a $k$-cube of $X_n(G)$ is determined uniquely by a set
$\{x_1,\ldots,x_k,u_1,\ldots,u_{n-k}\}$ of mutually disjoint simplices of $G$
where the $x_i$ are edges and the $u_i$ are vertices. The map $\Phi$ is
defined by sending the cube $\Cal C\{x_1,\ldots,x_k,u_1,\ldots,u_{n-k}\}$ onto the
cube in $\Cal T_{\Delta(G)}$ spanned by the edges 
$\overline x_1,\ldots,\overline x_k$ while respecting the orientations of all
edges.

\begin{figure}[ht!] 
\begin{center}
\begin{picture}(0,0)%
\includegraphics{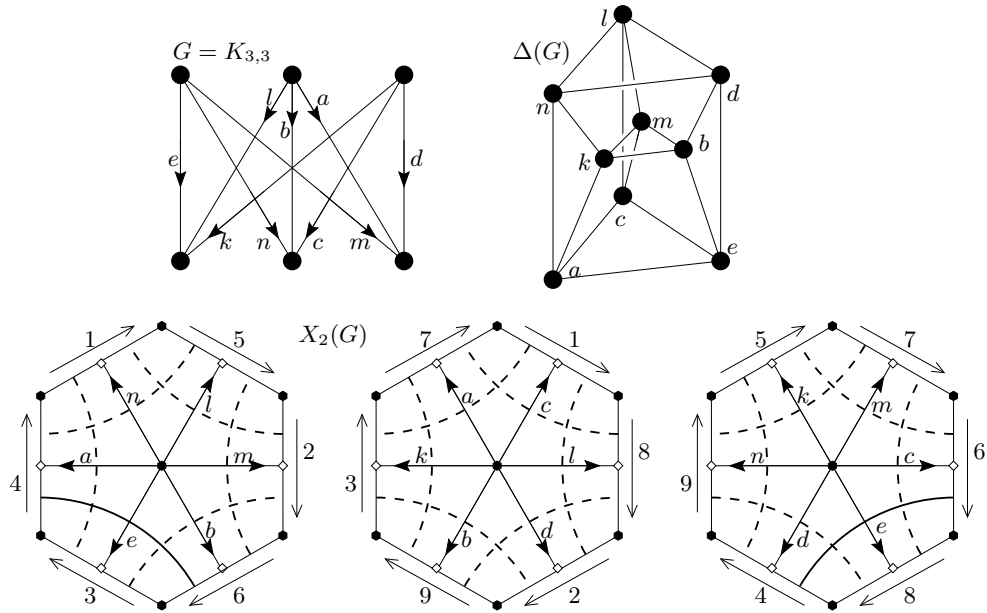}%
\end{picture}%
\setlength{\unitlength}{1539sp}%
\begingroup\makeatletter\ifx\SetFigFont\undefined%
\gdef\SetFigFont#1#2#3#4#5{%
  \reset@font\fontsize{#1}{#2pt}%
  \fontfamily{#3}\fontseries{#4}\fontshape{#5}%
  \selectfont}%
\fi\endgroup%
\begin{picture}(15525,9777)(226,-9298)
\put(2776,-2086){\makebox(0,0)[lb]{\smash{\SetFigFont{9}{10.8}{\rmdefault}{\mddefault}{\updefault}$e$}}}
\put(6676,-2086){\makebox(0,0)[lb]{\smash{\SetFigFont{9}{10.8}{\rmdefault}{\mddefault}{\updefault}$d$}}}
\put(4351,-1111){\makebox(0,0)[lb]{\smash{\SetFigFont{9}{10.8}{\rmdefault}{\mddefault}{\updefault}$l$}}}
\put(4576,-1636){\makebox(0,0)[lb]{\smash{\SetFigFont{9}{10.8}{\rmdefault}{\mddefault}{\updefault}$b$}}}
\put(3601,-3436){\makebox(0,0)[lb]{\smash{\SetFigFont{9}{10.8}{\rmdefault}{\mddefault}{\updefault}$k$}}}
\put(4201,-3436){\makebox(0,0)[lb]{\smash{\SetFigFont{9}{10.8}{\rmdefault}{\mddefault}{\updefault}$n$}}}
\put(5101,-3436){\makebox(0,0)[lb]{\smash{\SetFigFont{9}{10.8}{\rmdefault}{\mddefault}{\updefault}$c$}}}
\put(5701,-3436){\makebox(0,0)[lb]{\smash{\SetFigFont{9}{10.8}{\rmdefault}{\mddefault}{\updefault}$m$}}}
\put(5176,-1111){\makebox(0,0)[lb]{\smash{\SetFigFont{9}{10.8}{\rmdefault}{\mddefault}{\updefault}$a$}}}
\put(2851,-361){\makebox(0,0)[lb]{\smash{\SetFigFont{9}{10.8}{\rmdefault}{\mddefault}{\updefault}$G=K_{3,3}$}}}
\put(8326,-361){\makebox(0,0)[lb]{\smash{\SetFigFont{9}{10.8}{\rmdefault}{\mddefault}{\updefault}$\Delta(G)$}}}
\put(11776,-3511){\makebox(0,0)[lb]{\smash{\SetFigFont{9}{10.8}{\rmdefault}{\mddefault}{\updefault}$e$}}}
\put(9976,-3061){\makebox(0,0)[lb]{\smash{\SetFigFont{9}{10.8}{\rmdefault}{\mddefault}{\updefault}$c$}}}
\put(11776,-1036){\makebox(0,0)[lb]{\smash{\SetFigFont{9}{10.8}{\rmdefault}{\mddefault}{\updefault}$d$}}}
\put(11326,-1861){\makebox(0,0)[lb]{\smash{\SetFigFont{9}{10.8}{\rmdefault}{\mddefault}{\updefault}$b$}}}
\put(9226,-3886){\makebox(0,0)[lb]{\smash{\SetFigFont{9}{10.8}{\rmdefault}{\mddefault}{\updefault}$a$}}}
\put(10576,-1486){\makebox(0,0)[lb]{\smash{\SetFigFont{9}{10.8}{\rmdefault}{\mddefault}{\updefault}$m$}}}
\put(8701,-1261){\makebox(0,0)[lb]{\smash{\SetFigFont{9}{10.8}{\rmdefault}{\mddefault}{\updefault}$n$}}}
\put(9751,164){\makebox(0,0)[lb]{\smash{\SetFigFont{9}{10.8}{\rmdefault}{\mddefault}{\updefault}$l$}}}
\put(9376,-2161){\makebox(0,0)[lb]{\smash{\SetFigFont{9}{10.8}{\rmdefault}{\mddefault}{\updefault}$k$}}}
\put(226,-7336){\makebox(0,0)[lb]{\smash{\SetFigFont{9}{10.8}{\rmdefault}{\mddefault}{\updefault}$4$}}}
\put(4951,-6811){\makebox(0,0)[lb]{\smash{\SetFigFont{9}{10.8}{\rmdefault}{\mddefault}{\updefault}$2$}}}
\put(1426,-9136){\makebox(0,0)[lb]{\smash{\SetFigFont{9}{10.8}{\rmdefault}{\mddefault}{\updefault}$3$}}}
\put(3826,-9136){\makebox(0,0)[lb]{\smash{\SetFigFont{9}{10.8}{\rmdefault}{\mddefault}{\updefault}$6$}}}
\put(1426,-5011){\makebox(0,0)[lb]{\smash{\SetFigFont{9}{10.8}{\rmdefault}{\mddefault}{\updefault}$1$}}}
\put(3826,-5011){\makebox(0,0)[lb]{\smash{\SetFigFont{9}{10.8}{\rmdefault}{\mddefault}{\updefault}$5$}}}
\put(10351,-6811){\makebox(0,0)[lb]{\smash{\SetFigFont{9}{10.8}{\rmdefault}{\mddefault}{\updefault}$8$}}}
\put(6826,-9136){\makebox(0,0)[lb]{\smash{\SetFigFont{9}{10.8}{\rmdefault}{\mddefault}{\updefault}$9$}}}
\put(9226,-9136){\makebox(0,0)[lb]{\smash{\SetFigFont{9}{10.8}{\rmdefault}{\mddefault}{\updefault}$2$}}}
\put(6826,-5011){\makebox(0,0)[lb]{\smash{\SetFigFont{9}{10.8}{\rmdefault}{\mddefault}{\updefault}$7$}}}
\put(9226,-5011){\makebox(0,0)[lb]{\smash{\SetFigFont{9}{10.8}{\rmdefault}{\mddefault}{\updefault}$1$}}}
\put(15751,-6811){\makebox(0,0)[lb]{\smash{\SetFigFont{9}{10.8}{\rmdefault}{\mddefault}{\updefault}$6$}}}
\put(12226,-9136){\makebox(0,0)[lb]{\smash{\SetFigFont{9}{10.8}{\rmdefault}{\mddefault}{\updefault}$4$}}}
\put(14626,-9136){\makebox(0,0)[lb]{\smash{\SetFigFont{9}{10.8}{\rmdefault}{\mddefault}{\updefault}$8$}}}
\put(12226,-5011){\makebox(0,0)[lb]{\smash{\SetFigFont{9}{10.8}{\rmdefault}{\mddefault}{\updefault}$5$}}}
\put(14626,-5011){\makebox(0,0)[lb]{\smash{\SetFigFont{9}{10.8}{\rmdefault}{\mddefault}{\updefault}$7$}}}
\put(5626,-7336){\makebox(0,0)[lb]{\smash{\SetFigFont{9}{10.8}{\rmdefault}{\mddefault}{\updefault}$3$}}}
\put(11026,-7336){\makebox(0,0)[lb]{\smash{\SetFigFont{9}{10.8}{\rmdefault}{\mddefault}{\updefault}$9$}}}
\put(2101,-8236){\makebox(0,0)[lb]{\smash{\SetFigFont{9}{10.8}{\rmdefault}{\mddefault}{\updefault}$e$}}}
\put(3376,-8086){\makebox(0,0)[lb]{\smash{\SetFigFont{9}{10.8}{\rmdefault}{\mddefault}{\updefault}$b$}}}
\put(1351,-6886){\makebox(0,0)[lb]{\smash{\SetFigFont{9}{10.8}{\rmdefault}{\mddefault}{\updefault}$a$}}}
\put(7501,-8236){\makebox(0,0)[lb]{\smash{\SetFigFont{9}{10.8}{\rmdefault}{\mddefault}{\updefault}$b$}}}
\put(7501,-5911){\makebox(0,0)[lb]{\smash{\SetFigFont{9}{10.8}{\rmdefault}{\mddefault}{\updefault}$a$}}}
\put(8776,-8086){\makebox(0,0)[lb]{\smash{\SetFigFont{9}{10.8}{\rmdefault}{\mddefault}{\updefault}$d$}}}
\put(14176,-8011){\makebox(0,0)[lb]{\smash{\SetFigFont{9}{10.8}{\rmdefault}{\mddefault}{\updefault}$e$}}}
\put(14626,-6886){\makebox(0,0)[lb]{\smash{\SetFigFont{9}{10.8}{\rmdefault}{\mddefault}{\updefault}$c$}}}
\put(12901,-8236){\makebox(0,0)[lb]{\smash{\SetFigFont{9}{10.8}{\rmdefault}{\mddefault}{\updefault}$d$}}}
\put(6751,-6886){\makebox(0,0)[lb]{\smash{\SetFigFont{9}{10.8}{\rmdefault}{\mddefault}{\updefault}$k$}}}
\put(9226,-6886){\makebox(0,0)[lb]{\smash{\SetFigFont{9}{10.8}{\rmdefault}{\mddefault}{\updefault}$l$}}}
\put(12901,-5911){\makebox(0,0)[lb]{\smash{\SetFigFont{9}{10.8}{\rmdefault}{\mddefault}{\updefault}$k$}}}
\put(2101,-5911){\makebox(0,0)[lb]{\smash{\SetFigFont{9}{10.8}{\rmdefault}{\mddefault}{\updefault}$n$}}}
\put(12151,-6886){\makebox(0,0)[lb]{\smash{\SetFigFont{9}{10.8}{\rmdefault}{\mddefault}{\updefault}$n$}}}
\put(3826,-6886){\makebox(0,0)[lb]{\smash{\SetFigFont{9}{10.8}{\rmdefault}{\mddefault}{\updefault}$m$}}}
\put(14101,-6061){\makebox(0,0)[lb]{\smash{\SetFigFont{9}{10.8}{\rmdefault}{\mddefault}{\updefault}$m$}}}
\put(3376,-6061){\makebox(0,0)[lb]{\smash{\SetFigFont{9}{10.8}{\rmdefault}{\mddefault}{\updefault}$l$}}}
\put(8776,-6061){\makebox(0,0)[lb]{\smash{\SetFigFont{9}{10.8}{\rmdefault}{\mddefault}{\updefault}$c$}}}
\put(4876,-4936){\makebox(0,0)[lb]{\smash{\SetFigFont{9}{10.8}{\rmdefault}{\mddefault}{\updefault}$X_2(G)$}}}
\end{picture}
\end{center}
\caption{The reduced 2-point configuration space $X_2(K_{3,3})$ of the
$K_{3,3}$-graph is a nonorientable surface of Euler characteristic $-3$,
with a cubing by 18 squares (here grouped into three hexagons, which are
glued according to the arrows labelled 1 to 9). Its fundamental group
embeds in the right-angled Artin group $A_{\Delta(G)}$.}\label{F:chi-3}
\end{figure}

\begin{thm}
Let $G$ be any finite graph, and $n\in\N$. Then the map $\Phi \co X_n(G)$
$\to \Cal T_{\Delta(G)}$ described above is a local isometry and hence
$\pi_1$-injective. Thus the reduced graph braid group $RB_n(G)$ embeds 
in the right-angled Artin group $A_{\Delta(G)}$.  
\end{thm}

\begin{proof}
For simplicity write $X=X_n(G)$. Since we know that $\Cal T_{\Delta(G)}$
is locally CAT(0) it suffices, by Theorem \ref{GromovThm}, to show that
the link of each vertex of $X$ is mapped injectively onto a full
subcomplex of $L_{\Delta(G)}$. A vertex of $X$ is just a set $U$ of $n$
distinct vertices of $G$. Each vertex in the link $Lk(U,X)$ is uniquely
determined by an edge $e$ in $G$, precisely one of whose vertices, $u$
say, lies in $U$. This link vertex shall be written $e^+$ or $e^-$
depending on whether $e$ happens to be oriented away from or into the
graph vertex $u$.  Now the map induced by $\Phi$ sends a vertex $e^\varep$
in $Lk(U,X)$ to the vertex $\ov e^\varep$ in $L_\Delta$ (where
$\varep\in\{ +,-\}$) and so is clearly injective on vertices. Moreover one
easily observes that $e_1^{\varep_1},\ldots,e_n^{\varep_n}$ span a simplex in
$Lk(U,X)$ if and only if $e_1,\ldots,e_n$ are mutually disjoint, if and only
if $\ov e_1,\ldots,\ov e_n$ span a complete graph in $\Delta(G)$, if and only
if $\ov e_1^{\varep_1},\ldots,\ov e_n^{\varep_n}$ span a simplex in
$L_{\Delta(G)}$.  Thus the image of $Lk(U,X)$ is indeed a full subcomplex
of $L_{\Delta(G)}$. 
\end{proof}

We illustrate this Theorem with two examples
which will be used in the next Section.
The first example is the case where $G$ is taken to be
the bi-partite graph $K_{3,3}$, and $n=2$.
The graph $\Delta(K_{3,3})$ and the complex $X_2(K_{3,3})$
which is in fact a squaring of the
closed nonorientable surface of Euler characteristic $-3$,
are shown in Figure \ref{F:chi-3}. 
The dashed lines dual to the squaring together with the labels $a,b,\ldots,m,n$ 
constitute a ``cellular dissection'' of the surface, as defined in
Section \ref{Sect4} below, and describe the embedding of
$RB_2(K_{3,3})$ (the fundamental group of the surface)
into $A_{\Delta(K_{3,3})}$.

The second example is the case where $G=K_5$, the complete 
graph on $5$ vertices, and again $n=2$.
We leave the reader to verify that in fact $X_2(K_5)$ is a squaring 
of the closed nonorientable surface of Euler
characteristic $-5$ with $10$ vertices all of which have valence $6$ 
(there are therefore $15$ faces and $30$ edges).
The graph $\Delta(K_5)$ is in fact the image of the $1$-skeleton of 
the boundary of the dodecahedron under the 
$2$-fold quotient defined by the antipodal map. 


\section{Surface groups in right-angled Artin groups}\label{Sect4}


In this section, we classify the surface groups which can be embedded in a
right-angled Artin group. Let $RP^2$ denote the real projective plane,
$Kl$ the Klein bottle, and $S_{-1}$ the non-orientable surface of Euler
characteristic $-1$.

\begin{thm}
Let $S$ be a surface. 
\begin{itemize}
\item[\rm(i)] If $S\neq RP^2$, $Kl$, or $S_{-1}$, then there exists a
right-angled Artin group $A$ and a (quasi-isometric) embedding of
$\pi_1(S)$ as a subgroup of $A$.
\item[\rm(ii)] If $S=RP^2$, $Kl$, or $S_{-1}$, then $\pi_1(S)$ does not embed
as a subgroup of any right-angled Artin group.
\end{itemize}
\end{thm}

\begin{proof}
Note that the torus group $\Z\times\Z$ is already a right-angled Artin
group, so embeds in one.  Similary, the fundamental group of any
non-closed surface group is free of finite rank, so is already a
right-angled Artin group. In Proposition \ref{surfaces} we shall exhibit
quasi-isometric embeddings for each closed surface of Euler characteristic
less than or equal to $-2$. The fundamental group of the $2$-sphere is
trivial, so embeds trivially.

On the other hand, $\pi_1(RP^2)$ has order $2$, so does not embed in any
torsion free group. The Klein bottle group $\pi_1(Kl)$ has generalised
torsion and so is not biorderable (see for instance \cite{RolfW}).
Therefore $\pi_1(Kl)$ cannot embed in any right-angled Artin group since
right-angled Artin groups are known to be biorderable \cite{DT} (see also
Section \ref{Sect5}). Finally it follows from Theorem \ref{x2y2z2} below
(see Corollary \ref{noembed}) that the group $\pi_1(S_{-1})$ cannot be
embedded in any right-angled Artin group. \end{proof}

In order to construct embeddings of surfaces we shall appeal once again to
the techniques of Section \ref{Sect2}. In particular, every embedding is
realised geometrically by a geodesic embedding of one CAT(0) space in
another. Consequently, the group homomorphism obtained in each case is a
quasi-isometric embedding. We note that similar embeddings of surface
groups in right-angled Artin groups have been obtained in \cite{DSS}. The
emphasis in that paper is somewhat different however. They show that
``most'' right-angled Artin groups contain non-abelian surface groups but
are less concerned with realising a given surface group.

Let $S$ be compact surface. Let $c$ denote an orientation preserving
simple closed curve in the interior of $S$. A \emph{transverse
orientation} of $c$ is a choice of one or other of the two connected
components of the total space of the normal bundle to $c$ in $S$, and is
simply indicated by an arrow crossing the curve from one side to the
other. We extend this definition in the obvious way to allow $c$ to be a
boundary component of $S$, in which case a transverse orientation is
either ``into'' or ``out of'' the surface, and to arbitrary intervals
$(I,\partial I)$ properly embedded in $(S,\partial S)$. Transversely
oriented simple closed curves, properly embedded intervals and boundary
components shall be called collectively \emph{hypercurves} in $S$.

\begin{defn}\label{dissect}
Let $\Cal V$ be a finite set, $S$ a compact surface. A \emph{$\Cal
V$-dissection} of $S$ is a collection $\Cal H$ of hypercurves in $S$
satisfying
\begin{itemize}
\item[\rm(i)] any two hypercurves of $\Cal H$ are either disjoint or
intersect transversely in a discrete set of points, and no more than two
hypercurves meet at any one point;
\item[\rm(ii)] as well as being transversely oriented, each hypercurve in
$\Cal H$ is labelled with an element of $\Cal V$ in such a way that no two
hypercurves which intersect have the same label.
\end{itemize}
\end{defn}

Note that we are interested a priori in $\Cal V$-dissections of closed
surfaces (in which case all hypercurves are simple closed curves). However
the above definition will allow us to construct dissections of closed
surfaces by decomposing them into compact pieces with boundary and
dissecting these separately. Our interest in the above definition stems
from the following observation:

A $\Cal V$-dissection $\Cal H$ of a closed surface $S$, together with a
choice of basepoint $p\in S\setminus(\cup\Cal H)$, determines a
homomorphism of $\pi_1(S,p)$ to a right-angled Artin group as follows. Let
$\Delta(\Cal H)$ be the graph with vertex set $\Cal V$ and edge set $\Cal
E$ where $(u,v)\in\Cal E$ if and only if there exist a pair of
intersecting hypercurves in $\Cal H$ labelled $u$ and $v$ respectively.
Given a loop $\gamma$ at $p$ in $S$ which is in general position with
respect to $\Cal H$ (in particular we require that $\gamma$ crosses only
one hypercurve at a time), we define $w_\gamma$ to be the word in the
letters $\Cal V\cup\Cal V^{-1}$ obtained by reading off the labels of the
hypercurves traversed by $\gamma$ (the sign of each letter in $w_\gamma$
is determined by whether $\gamma$ traverses the hyperplane with or against
the transverse orientation). Then the homomorphism $\psi\co \pi_1(S,p)\to
A_{\Delta(\Cal H)}$ is well-defined by setting $\psi([\gamma])=w_\gamma$.

A necessary condition for $\psi$ to be injective is that $\Cal H$ be
cellular:

\begin{defn}
We say that a dissection $\Cal H$ of $S$ is \emph{cellular} if
$S\setminus(\cup\Cal H)$ is a disjoint union of open disks. 
\end{defn}

In this case, $\psi$ is induced by a cubical map. If the dissection $\Cal
H$ is cellular, then it determines in an obvious way a polyhedral
decomposition of $S$. We consider the squaring of $S$ which is
dual to this decomposition, and which we denote $X(S,\Cal H)$. 
Moreover, each edge of $X(S,\Cal H)$ is naturally 
labelled with orientation by an element of
$\Cal V$ according to the label and transverse orientation of the
hypercurve to which it is dual. Opposite edges of a square in $X(S,\Cal
H)$ are labelled similarly, and edges labelled $u$ and $v$ span a square
only if $(u,v)$ is an edge of the graph $\Delta(\Cal H)$. Thus there is a 
uniquely defined cubical map
$\Psi\co X(S,\Cal H)\to \Cal T_{\Delta(\Cal H)}$ which sends each edge in
$X(S,\Cal H)$ to the unique edge in $\Cal T_{\Delta(\Cal H)}$ with the
same label (respecting orientations). Evidently $\Psi$ induces $\psi$ at
the level of fundamental groups.

We consider now a closed surface $S$ of Euler characteristic 
$\chi(S)< -1$. We shall describe a cellular dissection $\Cal H$ of $S$
for which the map $\Psi\co X(S,\Cal H)\to \Cal T_{\Delta(\Cal H)}$ is
locally isometric. There are several cases to consider:

\begin{figure}[ht!] 
\centerline{
\begin{picture}(0,0)%
\includegraphics[width=.95\hsize]{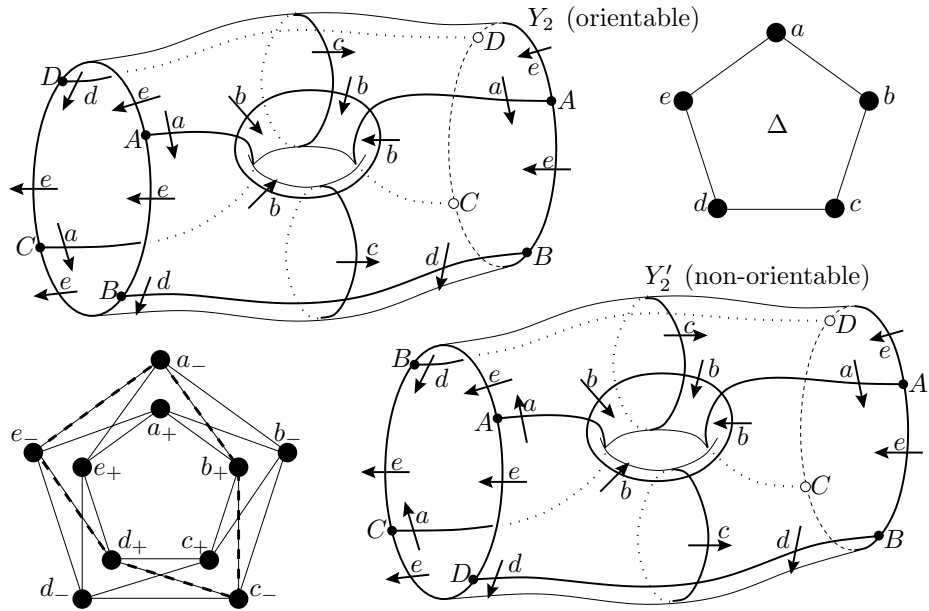}%
\end{picture}%
\setlength{\unitlength}{1612sp}%
\begingroup\makeatletter\ifx\SetFigFont\undefined%
\gdef\SetFigFont#1#2#3#4#5{%
  \reset@font\fontsize{#1}{#2pt}%
  \fontfamily{#3}\fontseries{#4}\fontshape{#5}%
  \selectfont}%
\fi\endgroup%
\begin{picture}(14091,9359)(193,-8760)
\put(12226, 89){\makebox(0,0)[lb]{\smash{\SetFigFont{10}{12.0}{\rmdefault}{\mddefault}{\updefault}$a$}}}
\put(13651,-961){\makebox(0,0)[lb]{\smash{\SetFigFont{10}{12.0}{\rmdefault}{\mddefault}{\updefault}$b$}}}
\put(10201,-961){\makebox(0,0)[lb]{\smash{\SetFigFont{10}{12.0}{\rmdefault}{\mddefault}{\updefault}$e$}}}
\put(10726,-2611){\makebox(0,0)[lb]{\smash{\SetFigFont{10}{12.0}{\rmdefault}{\mddefault}{\updefault}$d$}}}
\put(13126,-2611){\makebox(0,0)[lb]{\smash{\SetFigFont{10}{12.0}{\rmdefault}{\mddefault}{\updefault}$c$}}}
\put(11851,-1411){\makebox(0,0)[lb]{\smash{\SetFigFont{10}{12.0}{\rmdefault}{\mddefault}{\updefault}$\Delta$}}}
\put(7576,-811){\makebox(0,0)[lb]{\smash{\SetFigFont{10}{12.0}{\rmdefault}{\mddefault}{\updefault}$a$}}}
\put(4201,-2686){\makebox(0,0)[lb]{\smash{\SetFigFont{10}{12.0}{\rmdefault}{\mddefault}{\updefault}$b$}}}
\put(5176,-136){\makebox(0,0)[lb]{\smash{\SetFigFont{10}{12.0}{\rmdefault}{\mddefault}{\updefault}$c$}}}
\put(5701,-3286){\makebox(0,0)[lb]{\smash{\SetFigFont{10}{12.0}{\rmdefault}{\mddefault}{\updefault}$c$}}}
\put(1351,-961){\makebox(0,0)[lb]{\smash{\SetFigFont{10}{12.0}{\rmdefault}{\mddefault}{\updefault}$d$}}}
\put(2476,-2386){\makebox(0,0)[lb]{\smash{\SetFigFont{10}{12.0}{\rmdefault}{\mddefault}{\updefault}$e$}}}
\put(976,-3811){\makebox(0,0)[lb]{\smash{\SetFigFont{10}{12.0}{\rmdefault}{\mddefault}{\updefault}$e$}}}
\put(676,-2236){\makebox(0,0)[lb]{\smash{\SetFigFont{10}{12.0}{\rmdefault}{\mddefault}{\updefault}$e$}}}
\put(2176,-886){\makebox(0,0)[lb]{\smash{\SetFigFont{10}{12.0}{\rmdefault}{\mddefault}{\updefault}$e$}}}
\put(8401,-1936){\makebox(0,0)[lb]{\smash{\SetFigFont{10}{12.0}{\rmdefault}{\mddefault}{\updefault}$e$}}}
\put(8176,-511){\makebox(0,0)[lb]{\smash{\SetFigFont{10}{12.0}{\rmdefault}{\mddefault}{\updefault}$e$}}}
\put(1951,-1636){\makebox(0,0)[lb]{\smash{\SetFigFont{10}{12.0}{\rmdefault}{\mddefault}{\updefault}$A$}}}
\put(1576,-3961){\makebox(0,0)[lb]{\smash{\SetFigFont{10}{12.0}{\rmdefault}{\mddefault}{\updefault}$B$}}}
\put(301,-3286){\makebox(0,0)[lb]{\smash{\SetFigFont{10}{12.0}{\rmdefault}{\mddefault}{\updefault}$C$}}}
\put(676,-661){\makebox(0,0)[lb]{\smash{\SetFigFont{10}{12.0}{\rmdefault}{\mddefault}{\updefault}$D$}}}
\put(7501,-136){\makebox(0,0)[lb]{\smash{\SetFigFont{10}{12.0}{\rmdefault}{\mddefault}{\updefault}$D$}}}
\put(8626,-1036){\makebox(0,0)[lb]{\smash{\SetFigFont{10}{12.0}{\rmdefault}{\mddefault}{\updefault}$A$}}}
\put(8251,-3436){\makebox(0,0)[lb]{\smash{\SetFigFont{10}{12.0}{\rmdefault}{\mddefault}{\updefault}$B$}}}
\put(7126,-2611){\makebox(0,0)[lb]{\smash{\SetFigFont{10}{12.0}{\rmdefault}{\mddefault}{\updefault}$C$}}}
\put(1051,-3061){\makebox(0,0)[lb]{\smash{\SetFigFont{10}{12.0}{\rmdefault}{\mddefault}{\updefault}$a$}}}
\put(2476,-3811){\makebox(0,0)[lb]{\smash{\SetFigFont{10}{12.0}{\rmdefault}{\mddefault}{\updefault}$d$}}}
\put(6001,-1861){\makebox(0,0)[lb]{\smash{\SetFigFont{10}{12.0}{\rmdefault}{\mddefault}{\updefault}$b$}}}
\put(5551,-811){\makebox(0,0)[lb]{\smash{\SetFigFont{10}{12.0}{\rmdefault}{\mddefault}{\updefault}$b$}}}
\put(6601,-3436){\makebox(0,0)[lb]{\smash{\SetFigFont{10}{12.0}{\rmdefault}{\mddefault}{\updefault}$d$}}}
\put(2701,-1336){\makebox(0,0)[lb]{\smash{\SetFigFont{10}{12.0}{\rmdefault}{\mddefault}{\updefault}$a$}}}
\put(3676,-961){\makebox(0,0)[lb]{\smash{\SetFigFont{10}{12.0}{\rmdefault}{\mddefault}{\updefault}$b$}}}
\put(8176,239){\makebox(0,0)[lb]{\smash{\SetFigFont{10}{12.0}{\rmdefault}{\mddefault}{\updefault}$Y_2$ (orientable)}}}
\put(1501,-6661){\makebox(0,0)[lb]{\smash{\SetFigFont{10}{12.0}{\rmdefault}{\mddefault}{\updefault}$e_+$}}}
\put(1876,-7786){\makebox(0,0)[lb]{\smash{\SetFigFont{10}{12.0}{\rmdefault}{\mddefault}{\updefault}$d_+$}}}
\put(2851,-7786){\makebox(0,0)[lb]{\smash{\SetFigFont{10}{12.0}{\rmdefault}{\mddefault}{\updefault}$c_+$}}}
\put(3151,-6661){\makebox(0,0)[lb]{\smash{\SetFigFont{10}{12.0}{\rmdefault}{\mddefault}{\updefault}$b_+$}}}
\put(2326,-6061){\makebox(0,0)[lb]{\smash{\SetFigFont{10}{12.0}{\rmdefault}{\mddefault}{\updefault}$a_+$}}}
\put(4276,-6136){\makebox(0,0)[lb]{\smash{\SetFigFont{10}{12.0}{\rmdefault}{\mddefault}{\updefault}$b_-$}}}
\put(3901,-8536){\makebox(0,0)[lb]{\smash{\SetFigFont{10}{12.0}{\rmdefault}{\mddefault}{\updefault}$c_-$}}}
\put(676,-8536){\makebox(0,0)[lb]{\smash{\SetFigFont{10}{12.0}{\rmdefault}{\mddefault}{\updefault}$d_-$}}}
\put(226,-6136){\makebox(0,0)[lb]{\smash{\SetFigFont{10}{12.0}{\rmdefault}{\mddefault}{\updefault}$e_-$}}}
\put(2776,-5011){\makebox(0,0)[lb]{\smash{\SetFigFont{10}{12.0}{\rmdefault}{\mddefault}{\updefault}$a_-$}}}
\put(9976,-3736){\makebox(0,0)[lb]{\smash{\SetFigFont{10}{12.0}{\rmdefault}{\mddefault}{\updefault}$Y'_2$ (non-orientable)}}}
\put(12976,-5161){\makebox(0,0)[lb]{\smash{\SetFigFont{10}{12.0}{\rmdefault}{\mddefault}{\updefault}$a$}}}
\put(9601,-7036){\makebox(0,0)[lb]{\smash{\SetFigFont{10}{12.0}{\rmdefault}{\mddefault}{\updefault}$b$}}}
\put(10576,-4486){\makebox(0,0)[lb]{\smash{\SetFigFont{10}{12.0}{\rmdefault}{\mddefault}{\updefault}$c$}}}
\put(11101,-7636){\makebox(0,0)[lb]{\smash{\SetFigFont{10}{12.0}{\rmdefault}{\mddefault}{\updefault}$c$}}}
\put(6751,-5311){\makebox(0,0)[lb]{\smash{\SetFigFont{10}{12.0}{\rmdefault}{\mddefault}{\updefault}$d$}}}
\put(7876,-6736){\makebox(0,0)[lb]{\smash{\SetFigFont{10}{12.0}{\rmdefault}{\mddefault}{\updefault}$e$}}}
\put(6376,-8161){\makebox(0,0)[lb]{\smash{\SetFigFont{10}{12.0}{\rmdefault}{\mddefault}{\updefault}$e$}}}
\put(6076,-6586){\makebox(0,0)[lb]{\smash{\SetFigFont{10}{12.0}{\rmdefault}{\mddefault}{\updefault}$e$}}}
\put(7576,-5236){\makebox(0,0)[lb]{\smash{\SetFigFont{10}{12.0}{\rmdefault}{\mddefault}{\updefault}$e$}}}
\put(13801,-6286){\makebox(0,0)[lb]{\smash{\SetFigFont{10}{12.0}{\rmdefault}{\mddefault}{\updefault}$e$}}}
\put(13576,-4861){\makebox(0,0)[lb]{\smash{\SetFigFont{10}{12.0}{\rmdefault}{\mddefault}{\updefault}$e$}}}
\put(7351,-5986){\makebox(0,0)[lb]{\smash{\SetFigFont{10}{12.0}{\rmdefault}{\mddefault}{\updefault}$A$}}}
\put(6976,-8311){\makebox(0,0)[lb]{\smash{\SetFigFont{10}{12.0}{\rmdefault}{\mddefault}{\updefault}$D$}}}
\put(5701,-7636){\makebox(0,0)[lb]{\smash{\SetFigFont{10}{12.0}{\rmdefault}{\mddefault}{\updefault}$C$}}}
\put(6076,-5011){\makebox(0,0)[lb]{\smash{\SetFigFont{10}{12.0}{\rmdefault}{\mddefault}{\updefault}$B$}}}
\put(12901,-4486){\makebox(0,0)[lb]{\smash{\SetFigFont{10}{12.0}{\rmdefault}{\mddefault}{\updefault}$D$}}}
\put(14026,-5386){\makebox(0,0)[lb]{\smash{\SetFigFont{10}{12.0}{\rmdefault}{\mddefault}{\updefault}$A$}}}
\put(13651,-7786){\makebox(0,0)[lb]{\smash{\SetFigFont{10}{12.0}{\rmdefault}{\mddefault}{\updefault}$B$}}}
\put(12526,-6961){\makebox(0,0)[lb]{\smash{\SetFigFont{10}{12.0}{\rmdefault}{\mddefault}{\updefault}$C$}}}
\put(6451,-7411){\makebox(0,0)[lb]{\smash{\SetFigFont{10}{12.0}{\rmdefault}{\mddefault}{\updefault}$a$}}}
\put(7876,-8161){\makebox(0,0)[lb]{\smash{\SetFigFont{10}{12.0}{\rmdefault}{\mddefault}{\updefault}$d$}}}
\put(11401,-6211){\makebox(0,0)[lb]{\smash{\SetFigFont{10}{12.0}{\rmdefault}{\mddefault}{\updefault}$b$}}}
\put(10951,-5161){\makebox(0,0)[lb]{\smash{\SetFigFont{10}{12.0}{\rmdefault}{\mddefault}{\updefault}$b$}}}
\put(12001,-7786){\makebox(0,0)[lb]{\smash{\SetFigFont{10}{12.0}{\rmdefault}{\mddefault}{\updefault}$d$}}}
\put(8101,-5686){\makebox(0,0)[lb]{\smash{\SetFigFont{10}{12.0}{\rmdefault}{\mddefault}{\updefault}$a$}}}
\put(9076,-5311){\makebox(0,0)[lb]{\smash{\SetFigFont{10}{12.0}{\rmdefault}{\mddefault}{\updefault}$b$}}}
\end{picture}
}
\caption{All surfaces $S$ with even negative Euler characteristic can be
obtained by glueing copies of $Y_2$ and $Y'_2$. The cubical decomposition
$X(S,\Cal H)$ of $S$ contains one vertex in each component of 
$S\setminus(\cup\Cal H)$, and its link is a 5-cycle. The groups $\pi_1(S)$ 
all embed in $A_\Delta$, where $\Delta$ is as shown.
The link $L_\Delta$  of the vertex of $\Cal T_\Delta$ is a graph with
10 vertices and 20 edges, and contains the image of the link of any
vertex of $X(S,\Cal H)$ as a full subgraph. } \label{F:0mod2}
\end{figure}

\paragraph{Case $\chi(S)\equiv 0$ mod $2$}
 Let $\Cal V$ be the set $\{ a,b,c,d,e\}$ and consider the $\Cal
V$-dissected surfaces $Y_2$ and $Y_2'$ illustrated in figure
\ref{F:0mod2}. Note that these dissections are cellular. Also
 $\chi(Y_2)=\chi(Y_2')=-2$.  Identification of the two boundary components
of $Y_2$ while respecting the labelling of the vertices $A,B,C,D$ as shown
yields a $\Cal V$-dissection of the closed orientable surface of Euler
characteristic $-2$. A similar operation on $Y_2'$ yields cellular $\Cal
V$-dissection of the closed non-orientable surface of Euler characteristic
$-2$.  The graph $\Delta$ associated to each of these dissections is the
cycle of length $5$ with vertices labelled as in Figure \ref{F:0mod2}. It
is easily checked that the cubical map $\Psi$ sending the 
dual squaring into $\Cal T_\Delta$ satisfies, in each case, 
the conditions of Theorem \ref{GromovThm} (2), and so is a local isometry. 
(The link of each vertex in the squared surface is a $5$-cycle with
vertices labelled $a^\pm,b^\pm,c^\pm,d^\pm,e^\pm$ in that order around the
cycle).
 
More generally, we can glue $k$ copies of $Y_2$ end to end in a cycle
(while always respecting the labelling of vertices $A,B,C,D$) to construct
a $\Cal V$-dissection $\Cal H$ of the orientable surface $S$ of Euler
characteristic $-2k$ for any $k\geq 1$. By the argument just given, the
cubical map $\Psi\co X(S,\Cal H)\to \Cal T_{\Delta(\Cal H)}$ is locally
isometric.  To obtain a similar result for the non-orientable surface of
Euler characteristic $-2k$ we simply replace one instance of $Y_2$ with a
copy of $Y_2'$ in the above construction.

\begin{figure}[ht!] 
\centerline{
\begin{picture}(0,0)%
\includegraphics{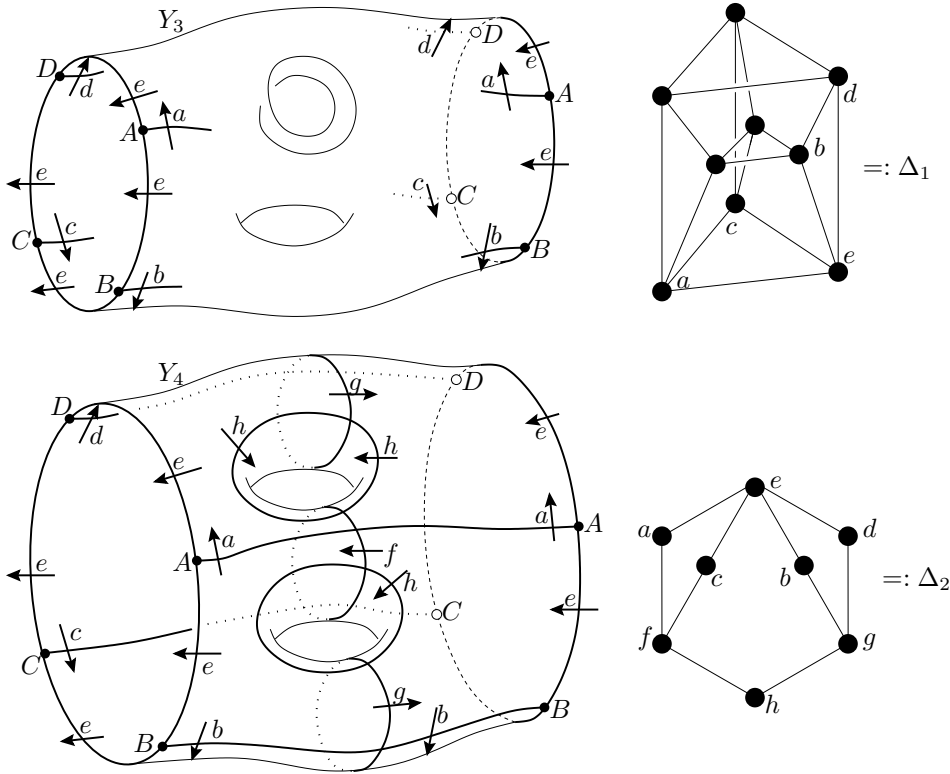}%
\end{picture}%
\setlength{\unitlength}{1618sp}%
\begingroup\makeatletter\ifx\SetFigFont\undefined%
\gdef\SetFigFont#1#2#3#4#5{%
  \reset@font\fontsize{#1}{#2pt}%
  \fontfamily{#3}\fontseries{#4}\fontshape{#5}%
  \selectfont}%
\fi\endgroup%
\begin{picture}(14820,11807)(193,-11344)
\put(1351,-961){\makebox(0,0)[lb]{\smash{\SetFigFont{10}{12.0}{\rmdefault}{\mddefault}{\updefault}$d$}}}
\put(2476,-2386){\makebox(0,0)[lb]{\smash{\SetFigFont{10}{12.0}{\rmdefault}{\mddefault}{\updefault}$e$}}}
\put(976,-3811){\makebox(0,0)[lb]{\smash{\SetFigFont{10}{12.0}{\rmdefault}{\mddefault}{\updefault}$e$}}}
\put(676,-2236){\makebox(0,0)[lb]{\smash{\SetFigFont{10}{12.0}{\rmdefault}{\mddefault}{\updefault}$e$}}}
\put(2176,-886){\makebox(0,0)[lb]{\smash{\SetFigFont{10}{12.0}{\rmdefault}{\mddefault}{\updefault}$e$}}}
\put(8401,-1936){\makebox(0,0)[lb]{\smash{\SetFigFont{10}{12.0}{\rmdefault}{\mddefault}{\updefault}$e$}}}
\put(8176,-511){\makebox(0,0)[lb]{\smash{\SetFigFont{10}{12.0}{\rmdefault}{\mddefault}{\updefault}$e$}}}
\put(1951,-1636){\makebox(0,0)[lb]{\smash{\SetFigFont{10}{12.0}{\rmdefault}{\mddefault}{\updefault}$A$}}}
\put(1576,-3961){\makebox(0,0)[lb]{\smash{\SetFigFont{10}{12.0}{\rmdefault}{\mddefault}{\updefault}$B$}}}
\put(301,-3286){\makebox(0,0)[lb]{\smash{\SetFigFont{10}{12.0}{\rmdefault}{\mddefault}{\updefault}$C$}}}
\put(676,-661){\makebox(0,0)[lb]{\smash{\SetFigFont{10}{12.0}{\rmdefault}{\mddefault}{\updefault}$D$}}}
\put(7501,-136){\makebox(0,0)[lb]{\smash{\SetFigFont{10}{12.0}{\rmdefault}{\mddefault}{\updefault}$D$}}}
\put(8626,-1036){\makebox(0,0)[lb]{\smash{\SetFigFont{10}{12.0}{\rmdefault}{\mddefault}{\updefault}$A$}}}
\put(8251,-3436){\makebox(0,0)[lb]{\smash{\SetFigFont{10}{12.0}{\rmdefault}{\mddefault}{\updefault}$B$}}}
\put(7126,-2611){\makebox(0,0)[lb]{\smash{\SetFigFont{10}{12.0}{\rmdefault}{\mddefault}{\updefault}$C$}}}
\put(2476,-3811){\makebox(0,0)[lb]{\smash{\SetFigFont{10}{12.0}{\rmdefault}{\mddefault}{\updefault}$b$}}}
\put(6451,-2386){\makebox(0,0)[lb]{\smash{\SetFigFont{10}{12.0}{\rmdefault}{\mddefault}{\updefault}$c$}}}
\put(7651,-3211){\makebox(0,0)[lb]{\smash{\SetFigFont{10}{12.0}{\rmdefault}{\mddefault}{\updefault}$b$}}}
\put(6526,-286){\makebox(0,0)[lb]{\smash{\SetFigFont{10}{12.0}{\rmdefault}{\mddefault}{\updefault}$d$}}}
\put(3376,-10786){\makebox(0,0)[lb]{\smash{\SetFigFont{10}{12.0}{\rmdefault}{\mddefault}{\updefault}$b$}}}
\put(901,-5836){\makebox(0,0)[lb]{\smash{\SetFigFont{10}{12.0}{\rmdefault}{\mddefault}{\updefault}$D$}}}
\put(1501,-6286){\makebox(0,0)[lb]{\smash{\SetFigFont{10}{12.0}{\rmdefault}{\mddefault}{\updefault}$d$}}}
\put(8326,-7486){\makebox(0,0)[lb]{\smash{\SetFigFont{10}{12.0}{\rmdefault}{\mddefault}{\updefault}$a$}}}
\put(6826,-10561){\makebox(0,0)[lb]{\smash{\SetFigFont{10}{12.0}{\rmdefault}{\mddefault}{\updefault}$b$}}}
\put(2776,-8236){\makebox(0,0)[lb]{\smash{\SetFigFont{10}{12.0}{\rmdefault}{\mddefault}{\updefault}$A$}}}
\put(2176,-11011){\makebox(0,0)[lb]{\smash{\SetFigFont{10}{12.0}{\rmdefault}{\mddefault}{\updefault}$B$}}}
\put(451,-9811){\makebox(0,0)[lb]{\smash{\SetFigFont{10}{12.0}{\rmdefault}{\mddefault}{\updefault}$C$}}}
\put(8551,-10486){\makebox(0,0)[lb]{\smash{\SetFigFont{10}{12.0}{\rmdefault}{\mddefault}{\updefault}$B$}}}
\put(6901,-8986){\makebox(0,0)[lb]{\smash{\SetFigFont{10}{12.0}{\rmdefault}{\mddefault}{\updefault}$C$}}}
\put(9076,-7636){\makebox(0,0)[lb]{\smash{\SetFigFont{10}{12.0}{\rmdefault}{\mddefault}{\updefault}$A$}}}
\put(7201,-5461){\makebox(0,0)[lb]{\smash{\SetFigFont{10}{12.0}{\rmdefault}{\mddefault}{\updefault}$D$}}}
\put(2776,-6661){\makebox(0,0)[lb]{\smash{\SetFigFont{10}{12.0}{\rmdefault}{\mddefault}{\updefault}$e$}}}
\put(676,-8236){\makebox(0,0)[lb]{\smash{\SetFigFont{10}{12.0}{\rmdefault}{\mddefault}{\updefault}$e$}}}
\put(1351,-10711){\makebox(0,0)[lb]{\smash{\SetFigFont{10}{12.0}{\rmdefault}{\mddefault}{\updefault}$e$}}}
\put(3226,-9811){\makebox(0,0)[lb]{\smash{\SetFigFont{10}{12.0}{\rmdefault}{\mddefault}{\updefault}$e$}}}
\put(9901,-9361){\makebox(0,0)[lb]{\smash{\SetFigFont{10}{12.0}{\rmdefault}{\mddefault}{\updefault}$f$}}}
\put(9901,-7711){\makebox(0,0)[lb]{\smash{\SetFigFont{10}{12.0}{\rmdefault}{\mddefault}{\updefault}$a$}}}
\put(11026,-8461){\makebox(0,0)[lb]{\smash{\SetFigFont{10}{12.0}{\rmdefault}{\mddefault}{\updefault}$c$}}}
\put(11926,-6961){\makebox(0,0)[lb]{\smash{\SetFigFont{10}{12.0}{\rmdefault}{\mddefault}{\updefault}$e$}}}
\put(13351,-7711){\makebox(0,0)[lb]{\smash{\SetFigFont{10}{12.0}{\rmdefault}{\mddefault}{\updefault}$d$}}}
\put(13351,-9361){\makebox(0,0)[lb]{\smash{\SetFigFont{10}{12.0}{\rmdefault}{\mddefault}{\updefault}$g$}}}
\put(12076,-8461){\makebox(0,0)[lb]{\smash{\SetFigFont{10}{12.0}{\rmdefault}{\mddefault}{\updefault}$b$}}}
\put(6001,-8086){\makebox(0,0)[lb]{\smash{\SetFigFont{10}{12.0}{\rmdefault}{\mddefault}{\updefault}$f$}}}
\put(6301,-8536){\makebox(0,0)[lb]{\smash{\SetFigFont{10}{12.0}{\rmdefault}{\mddefault}{\updefault}$h$}}}
\put(11851,-10411){\makebox(0,0)[lb]{\smash{\SetFigFont{10}{12.0}{\rmdefault}{\mddefault}{\updefault}$h$}}}
\put(6001,-6436){\makebox(0,0)[lb]{\smash{\SetFigFont{10}{12.0}{\rmdefault}{\mddefault}{\updefault}$h$}}}
\put(5476,-5461){\makebox(0,0)[lb]{\smash{\SetFigFont{10}{12.0}{\rmdefault}{\mddefault}{\updefault}$g$}}}
\put(6151,-10186){\makebox(0,0)[lb]{\smash{\SetFigFont{10}{12.0}{\rmdefault}{\mddefault}{\updefault}$g$}}}
\put(3676,-6136){\makebox(0,0)[lb]{\smash{\SetFigFont{10}{12.0}{\rmdefault}{\mddefault}{\updefault}$h$}}}
\put(8776,-8761){\makebox(0,0)[lb]{\smash{\SetFigFont{10}{12.0}{\rmdefault}{\mddefault}{\updefault}$e$}}}
\put(8326,-6211){\makebox(0,0)[lb]{\smash{\SetFigFont{10}{12.0}{\rmdefault}{\mddefault}{\updefault}$e$}}}
\put(13351,-2161){\makebox(0,0)[lb]{\smash{\SetFigFont{10}{12.0}{\rmdefault}{\mddefault}{\updefault}$=: \Delta_1$}}}
\put(13651,-8461){\makebox(0,0)[lb]{\smash{\SetFigFont{10}{12.0}{\rmdefault}{\mddefault}{\updefault}$=: \Delta_2$}}}
\put(2551, 89){\makebox(0,0)[lb]{\smash{\SetFigFont{10}{12.0}{\rmdefault}{\mddefault}{\updefault}$Y_3$}}}
\put(2551,-5311){\makebox(0,0)[lb]{\smash{\SetFigFont{10}{12.0}{\rmdefault}{\mddefault}{\updefault}$Y_4$}}}
\put(13051,-3511){\makebox(0,0)[lb]{\smash{\SetFigFont{10}{12.0}{\rmdefault}{\mddefault}{\updefault}$e$}}}
\put(11251,-3061){\makebox(0,0)[lb]{\smash{\SetFigFont{10}{12.0}{\rmdefault}{\mddefault}{\updefault}$c$}}}
\put(13051,-1036){\makebox(0,0)[lb]{\smash{\SetFigFont{10}{12.0}{\rmdefault}{\mddefault}{\updefault}$d$}}}
\put(12601,-1861){\makebox(0,0)[lb]{\smash{\SetFigFont{10}{12.0}{\rmdefault}{\mddefault}{\updefault}$b$}}}
\put(10501,-3886){\makebox(0,0)[lb]{\smash{\SetFigFont{10}{12.0}{\rmdefault}{\mddefault}{\updefault}$a$}}}
\put(2776,-1336){\makebox(0,0)[lb]{\smash{\SetFigFont{10}{12.0}{\rmdefault}{\mddefault}{\updefault}$a$}}}
\put(3526,-7861){\makebox(0,0)[lb]{\smash{\SetFigFont{10}{12.0}{\rmdefault}{\mddefault}{\updefault}$a$}}}
\put(1201,-9286){\makebox(0,0)[lb]{\smash{\SetFigFont{10}{12.0}{\rmdefault}{\mddefault}{\updefault}$c$}}}
\put(1126,-3061){\makebox(0,0)[lb]{\smash{\SetFigFont{10}{12.0}{\rmdefault}{\mddefault}{\updefault}$c$}}}
\put(7501,-811){\makebox(0,0)[lb]{\smash{\SetFigFont{10}{12.0}{\rmdefault}{\mddefault}{\updefault}$a$}}}
\end{picture}
}
\caption{All surfaces $S$ with $\chi(S) \equiv -3$ mod $4$, where
$\chi(S)\leqslant -3$, and with $\chi(S) \equiv -1$ mod $4$, where
$\chi(S)\leqslant -9$, can be obtained by glueing copies of $Y_3$ and
$Y_4$. Their fundamental groups embed in the right-angled Artin group
$\Cal T_{\Delta(\Cal H)}$, where $\Delta(\Cal H)$ is the graph obtained by
identifying the vertices $a,b,c,d$ and $e$, and the edges connecting them,
in $\Delta_1$ and $\Delta_2$.}\label{F:3mod4} 
\end{figure}

\paragraph{Case $\chi(S)\equiv -3$ mod $4$}
Recall that the the cubical complex $X_2(K_{3,3})$ is simply the squaring of
the closed surface $S_{-3}$ of Euler characteristic $-3$ 
illustrated in Figure \ref{F:chi-3}.
Dual to this is a $\Cal V_1$-dissection $\Cal H_1$ of the closed surface (also
illustrated in Figure  \ref{F:chi-3}), where $\Cal V_1=\{a,b,c,d,e,k,l,m,n\}$, 
and where $\Delta(\Cal H_1)=\Delta(K_{3,3})$ is the graph $\Delta_1$ shown
in Figure \ref{F:3mod4}. Note that every hypercurve of $\Cal H_1$
has edge-length $4$ (as measured in the polygonal decomposition of $S_{-3}$
associated to this dissection). Also, there are no two hypercurves with the
same label.  Now cut $S_{-3}$ open along the hypercurve labelled $e$  to
obtain a surface with two boundary components whose vertices are labelled
$A,B,C,D$ as in Figure \ref{F:3mod4}. 
This surface, which we denote $Y_3$, carries a naturally
induced $\Cal V_1$-dissection. Consider also the $\Cal V_2$-dissected
surface $Y_4$ shown in Figure \ref{F:3mod4}, where $\Cal V_2=\{
a,b,c,d,e,f,g,h\}$ are the vertices of the associated graph $\Delta_2$
shown in the figure. Let $\Cal V=\Cal V_1\cup\Cal V_2$, where the labels
$a,b,c,d,e$ are the only elements common to $\Cal V_1$ and $\Cal V_2$.

Now gluing $k$ copies of $Y_4$ and one copy of $Y_3$ end to end in a cycle
(while respecting, of course, the transverse orientations on the boundary
components as well as the vertex labels $A,B,C,D$) results in a cellular
$\Cal V$-dissection $\Cal H$ of the surface $S$ of Euler characteristic
$-(3+4k)$, for $k\geq 0$.  It is again easily seen that the cubical map
$\Psi\co X(S,\Cal H)\to \Cal T_{\Delta(\Cal H)}$ is locally isometric.
Note that $\Delta(\Cal H)$ is the graph obtained by amalgamating
$\Delta_1$ and $\Delta_2$ as shown in Figure \ref{F:3mod4}.

\paragraph{Case $\chi(S)\equiv -1$ mod $4$}
We suppose $\chi(S)<-1$. The case $\chi(S)=-5$ is covered by the fact
that $X_2(K_5)$ is a squaring of this surface which embeds in $\Cal
T_{\Delta(K_5)}$ via a locally isometric cubical map, as in Section
\ref{Sect3}. The remaining cases $\chi(S)\leq -9$ may be treated exactly
as in the previous paragraph by gluing together $3$ copies of $Y_3$ and
$k$ copies of $Y_4$ to obtain a cellular dissection of the surface of
Euler characteristic $-(9+4k)$, for $k\geq 0$. \medskip

Combining the three preceding paragraphs gives:

\begin{prop}\label{surfaces}

For each closed surface $S$ other than $RP^2$, $Kl$, and $S_{-1}$, there
exists a squaring $X(S)$ of $S$, a graph $\Delta$ (without loops or
multiple edges), and a locally isometric cubical map 
$\Psi\co X(S)\to \Cal T_\Delta$, inducing an injection 
$\psi\co \pi_1(S)\to A_\Delta$ on the level of fundamental groups.
\noproof \end{prop}

We turn now to consider the surface $S_{-1}$ of Euler characteristic $-1$,
whose fundamental group admits the following presentation:
\[
\pi_1(S_{-1})=\< X,Y,Z\mid X^2Y^2=Z^2\>\,.
\]
The following Theorem generalises a well-known Theorem of Lyndon \cite{Lyn}
for free groups. See \cite{LSbook} for a more complete discussion
of the equation $x^2y^2=z^2$ over a free group.

\begin{thm}\label{x2y2z2}
Let $A$ denote an arbitrary right-angled Artin group. If elements
$x,y,z\in A$ satisfy the relation $x^2y^2=z^2$ in $A$ then $x$, $y$ and
$z$ are mutually commuting elements. \end{thm}

\begin{proof}
 We will use the fact that the relation $x^2y^2=z^2$ in $A$ implies the
existence of a homomorphism $\varphi\co \pi_1(S_{-1})\to A$ such that
$\varphi(X)=x$, $\varphi(Y)=y$ and $\varphi(Z)=z$. In the preceding
paragraphs we constructed homomorphisms from surface groups to
right-angled Artin groups by producing dissections of the surface. The key
observation at this point is that this method is completely general:

\paragraph{Claim A} \textsl{Let $S$ be a closed surface, and $A$ a
right-angled Artin group. Every homomorphism $\psi\co \pi_1(S)\to A$ is
defined by a dissection of $S$.}

\begin{proof}[Proof of Claim] 
 Consider a polygonal decomposition of the surface $S$ with one vertex and
one $2$-cell $D$, and orient the edges. This gives a $1$-relator
presentation $\Cal P$ of the closed surface group $\pi_1(S)$ where the
generators correspond to the labelled edges and the relator word can be
read around the boundary of the polygon $D$. (e.g: the presentation
already given above for $\pi_1(S_{-1})$ comes from description of $S_{-1}$
as a hexagon with face pairings).

Let $\Cal V$ denote the standard generating set of $A$ (the vertex set of
the defining graph $\Delta$). For each generator $U$ of $\Cal P$ choose a
word in the letters $\Cal V$ and their inverses which represents $\psi(U)$
(and which, by abuse of notation, we shall also denote by $U$).  Now
subdivide the boundary of the polygon $D$ so that its sides are labelled
now by the \emph{words} $X, Y, Z,\ldots,$ the new edges being oriented and
labelled by the letters of $\Cal V$. 

We now show how to construct a dissection of the surface $D$ which is 
\emph{consistent with the boundary labelling}: that is to say that
(i) the endpoints of those hypercurves which are properly embedded arcs
lie on the midpoints of edges of $\partial D$ and the labelling and
tranverse orientation of each such hypercurve agrees with the edge
labelling at each of its endpoints, and (ii) every edge
in $\partial D$ is at the endpoint of exactly one hypercurve.  

\begin{figure}[ht!] 
\begin{center}
\begin{picture}(0,0)%
\includegraphics{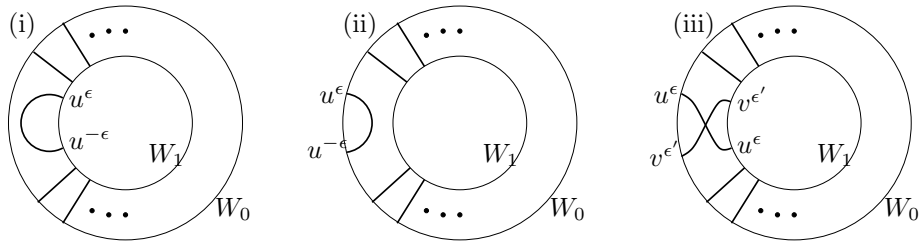}%
\end{picture}%
\setlength{\unitlength}{1381sp}%
\begingroup\makeatletter\ifx\SetFigFont\undefined%
\gdef\SetFigFont#1#2#3#4#5{%
  \reset@font\fontsize{#1}{#2pt}%
  \fontfamily{#3}\fontseries{#4}\fontshape{#5}%
  \selectfont}%
\fi\endgroup%
\begin{picture}(16216,4214)(293,-3668)
\put(3976,-3211){\makebox(0,0)[lb]{\smash{\SetFigFont{10}{12.0}{\rmdefault}{\mddefault}{\updefault}$W_0$}}}
\put(9976,-3211){\makebox(0,0)[lb]{\smash{\SetFigFont{10}{12.0}{\rmdefault}{\mddefault}{\updefault}$W_0$}}}
\put(15976,-3211){\makebox(0,0)[lb]{\smash{\SetFigFont{10}{12.0}{\rmdefault}{\mddefault}{\updefault}$W_0$}}}
\put(301, 89){\makebox(0,0)[lb]{\smash{\SetFigFont{10}{12.0}{\rmdefault}{\mddefault}{\updefault}(i)}}}
\put(1351,-2011){\makebox(0,0)[lb]{\smash{\SetFigFont{10}{12.0}{\rmdefault}{\mddefault}{\updefault}$u^{-\epsilon}$}}}
\put(2776,-2236){\makebox(0,0)[lb]{\smash{\SetFigFont{10}{12.0}{\rmdefault}{\mddefault}{\updefault}$W_1$}}}
\put(1351,-1261){\makebox(0,0)[lb]{\smash{\SetFigFont{10}{12.0}{\rmdefault}{\mddefault}{\updefault}$u^\epsilon$}}}
\put(6301, 89){\makebox(0,0)[lb]{\smash{\SetFigFont{10}{12.0}{\rmdefault}{\mddefault}{\updefault}(ii)}}}
\put(8776,-2236){\makebox(0,0)[lb]{\smash{\SetFigFont{10}{12.0}{\rmdefault}{\mddefault}{\updefault}$W_1$}}}
\put(5851,-1186){\makebox(0,0)[lb]{\smash{\SetFigFont{10}{12.0}{\rmdefault}{\mddefault}{\updefault}$u^\epsilon$}}}
\put(5626,-2161){\makebox(0,0)[lb]{\smash{\SetFigFont{10}{12.0}{\rmdefault}{\mddefault}{\updefault}$u^{-\epsilon}$}}}
\put(14776,-2236){\makebox(0,0)[lb]{\smash{\SetFigFont{10}{12.0}{\rmdefault}{\mddefault}{\updefault}$W_1$}}}
\put(11851,-1186){\makebox(0,0)[lb]{\smash{\SetFigFont{10}{12.0}{\rmdefault}{\mddefault}{\updefault}$u^\epsilon$}}}
\put(11776,-2311){\makebox(0,0)[lb]{\smash{\SetFigFont{10}{12.0}{\rmdefault}{\mddefault}{\updefault}$v^{\epsilon'}$}}}
\put(13351,-1336){\makebox(0,0)[lb]{\smash{\SetFigFont{10}{12.0}{\rmdefault}{\mddefault}{\updefault}$v^{\epsilon'}$}}}
\put(13351,-2086){\makebox(0,0)[lb]{\smash{\SetFigFont{10}{12.0}{\rmdefault}{\mddefault}{\updefault}$u^\epsilon$}}}
\put(12226, 89){\makebox(0,0)[lb]{\smash{\SetFigFont{10}{12.0}{\rmdefault}{\mddefault}{\updefault}(iii)}}}
\end{picture}
\end{center}
\caption{Extending the dissection of $D'$ across the annulus $A$ to a
dissection of $D=D'\cup A$. }\label{F:annuli}
\end{figure}

The boundary word $W_0$ which is read around the relator disk $D$
($W_0=X^2Y^2Z^{-2}$ in the case $S_{-1}$) represents the identity in $A$.
There is therefore a sequence $W_0,W_1,W_2,\ldots,W_n=1$ of words in $\Cal
V\cup\Cal V^{-1}$ such that $W_i$ is obtained from $W_{i-1}$ by either (i)
a trivial insertion of $uu^{-1}$ or $u^{-1}u$ (for $u\in\Cal V$), or (ii)
a trivial deletion, or (iii) a substitution of a subword $u^\epsilon
v^{\epsilon'}$ of $W_i$ ($\epsilon,\epsilon'\in\{\pm 1\}$)
 with the word $v^{\epsilon'}u^\epsilon$ in the case that $(u,v)$ is an
edge of the defining graph $\Delta$. Decompose $D$ as the union of a disk
$D'$ in the interior of $D$ and an annulus $A$, and label the boundary of
$D'$ by the word $W_1$. By induction on the number of steps of type
(i)--(iii) needed to transform a word into the trivial word, we may
suppose that $D'$ admits a dissection which is consistent with the
labelling $W_1$ of the boundary. This dissection extends across the
annulus $A$ as shown in Figure \ref{F:annuli} (i), (ii) or (iii),
according to whether the transformation from $W_0$ to $W_1$ is of type
(i), (ii), or (iii) respectively, to give the desired dissection of the
disk $D$.

Having shown that there exists a dissection of $D$ consistent with the 
boundary labels, it remains to observe that such a dissection induces
a dissection of the surface $S$ (after pairwise identification of the 
edges of the polygon $D$). Moreover, this dissection of $S$ clearly
induces the given homomorphism $\psi : \pi_1(S)\to A$, as required. 
\end{proof}

Note that the definition of the homomorphism $\psi$ by a dissection $\Cal
H$ of a closed surface $S$ actually depends on a choice of basepoint
$\ast$ in $S$. In the above argument, the basepoint is simply just the
vertex of the polygonal decomposition of $S$.  Allowing the basepoint to
move with respect to the dissection $\Cal H$ only changes the definition
of $\psi$ up to an inner automorphism of the image group $A$. We say that
homomorphisms $\psi$ and $\psi'$ from $\pi_1(S)$ to $A$ are
\emph{equivalent} if there exists $g\in A$ such that
$\psi'(w)=g\psi(w)g^{-1}$ for all $w\in\pi_1(S)$.

\paragraph{Claim B}
{\sl Let $\Cal H$ be a dissection of $S$ inducing the homomorphism
$\psi$. The following modifications of $\Cal H$ do not change the
equivalence class of $\psi$.
\begin{itemize} 
 \item[\rm(i)] Removing from $\Cal H$ any or all  hypercurves which are
homotopically trivial in $S$;
 \item[\rm(ii)] Changing $\Cal H$ to a different dissection $\Cal H'$ by a
homotopy of one or more hypercurves in $\Cal H$ without changing labels
or orientations and without introducing any intersection between
hypercurves which are already disjoint in $\Cal H$. 
\end{itemize} }

\begin{proof}[Proof of Claim]
 (i) Let $c$ be a homotopically trivial hypercurve of $\cal H$. Then $c$
bounds a disk $B$ in $S$. Up to equivalence of $\psi$ we are free to
choose the basepoint of $S$ to lie outside $B$. Now observe that any
element of $\pi_1(S)$ may be represented by a loop which is disjoint from
from $B$, so that the curve $c$ does not enter at all into the definition
of $\psi$ and may be removed.

(ii) We leave part (ii) as an exercise for the reader.
\end{proof}

We now return to the homomorphism $\varphi\co \pi_1(S_{-1})\to A$
introduced at the beginning of the proof. Theorem \ref{x2y2z2} is now
proved by showing that $\varphi(XYZ^{-1})=1$. For then we have relations
$xy=z$ and $x^2y^2=z^2$ in $A$, from which it easily follows that $xy=yx$,
and moreover that both $x$ and $y$ commute also with $z$. (Note that this
proof only depends on the equivalence class of $\varphi$).

Choose a simple loop $\gamma$ at the basepoint in $S_{-1}$ which
represents the element $XYZ^{-1}$. We make the following important
observation:

\paragraph{Claim C}
{\sl Any homotopically nontrivial orientation preserving simple closed
curve $c$ in $S_{-1}$ which intersects $\gamma$ transversely must
intersect $\gamma$ in exactly two points and must bound a M\"obius band in
$S_{-1}$.}

\begin{proof}[Proof of Claim]
We cut $S_{-1}$ open along the simple closed curve $\gamma$. Since $\gamma$
is an orientation reversing curve in $S_{-1}$, the surface $T$ resulting
from this operation is homeomorphic to a torus with an open
disk removed. Let $b$ denote the boundary of $T$. Then $S_{-1}$ is
obtained as the quotient of $T$ by the antipodal map on $b$. The curve
$c$, as viewed in $T$, consists of a collection of mutually disjoint arcs
$c_1,c_2,\ldots,c_k$ embedded in $T$ so that their endpoints lie on distinct
points of $b$.

Define the simple closed curves $\lambda=S^1\times\{ 1\}$, 
 $\mu=\{ 1\}\times S^1$ and $\delta=\{ (x,x): x\in S^1\}$ in the torus
$S^1\times S^1$. These curves are disjoint except at the common point
$(1,1)$. Suppose now that $T$ is obtained from $S^1\times S^1$ by deleting
an open neighbourhood of $(1,1)$, and let $l,m,d$ denote the arcs in $T$
obtained by intersecting $\lambda$, $\mu$ and $\delta$ respectively with
$T$.

Up to a homeomorphism of $T$ we may suppose that each of the mutually
disjoint arcs $c_1,\ldots,c_k$ is parallel to one of $l$, $m$, or $d$, and
that at least one of them is parallel to $l$. Re-label these arcs
$l_1,\ldots,l_r,d_{r+1},\ldots,d_s,m_{s+1},\ldots,m_k$ accordingly (where 
$1\leq r\leq s\leq k$). We label the endpoints of an arc $\alpha$ by
$\alpha^+$ and $\alpha^-$. Now observe that, up to re-numbering and
re-orienting of the arcs, their endpoints appear in the following cyclic
order around $b$:
\begin{align*}
l_1^+,l_2^+,\ldots,l_r^+,d_{r+1}^+,\ldots,d_s^+,m_{s+1}^+,\ldots,m_k^+, 
l_r^-,\ldots,l_2^-,l_1^-, &d_s^-,\ldots,\\&d_{r+1}^-,m_k^-,\ldots,m_{s+1}^-\,.
\end{align*}
Since $c$ is formed from these arcs by identifying pairs of vertices which
are antipodal in this cyclic order, the arcs $l_1$ and $l_r$ (or in the
case $r=1$, just the arc $l_1$ by itself) form a connected component of
$c$. But since $c$ is connected and orientation preserving we must have
$r=s=k=2$ and $c$ consists of just $l_1\cup l_2$. Moreover, since $l_1$
and $l_2$ are parallel in $T$, $c$ bounds a M\"obius band in $S_{-1}$.
\end{proof}

\proof[Proof that $\varphi(XYZ^{-1})=1$]
 We may clearly consider $\varphi$ up to equivalence. By Claim~A,
$\varphi$ may be defined by a dissection $\Cal H$ of $S_{-1}$. By
Claim~B(i) we may suppose that every hypercurve of $\Cal H$ is
homotopically non-trivial. We may suppose also that the loop $\gamma$
which represents $XYZ^{-1}$ is transverse to $\Cal H$. Let $c\in\Cal H$ be
a hypercurve which intersects $\gamma$ in at least one point. Then, by
Claim C, $c$ bounds a M\"obius band $M$ and intersects $\gamma$ in two
points.  By choosing an innermost such hypercurve $c$ we may suppose that
$M$ does not contain any hypercurve in its interior.  By a homotopy of $c$
satisfying the conditions of Claim B(ii) we may futhermore assume up to an
equivalence of $\varphi$ that the following hold:
 \begin{itemize}
 \item that the basepoint $\ast$ of $S_{-1}$ does not lie in
$M$,
 \item that no two hypercurves of $\Cal H$ intersect on the
interior of $M$, and
 \item that any hypercurve $c'$ of $\Cal H$ which intersects
$M$ does so essentially -- no component of $c'\cap M$ cobounds a disk
together with a subinterval of $c$.
 \end{itemize}
 As a consequence of the last two of these conditions, any loop $\alpha$
at $\ast$ in $S_{-1}$ is homotopic to one which does not intersect any
hypercurves while in the interior of $M$. It now follows that deleting $c$
from the dissection $\Cal H$ does not change the definition of $\varphi$
(up to equivalence).

By a straightforward induction we may now suppose that there are no
hypercurves in $\Cal H$ which intersect $\gamma$. That is to say,
$\varphi([\gamma])=1$ as required. This completes the proof of Theorem
\ref{x2y2z2}.\end{proof}

\begin{cor}\label{noembed}
 The fundamental group of the surface $S_{-1}$ of Euler characteristic
$-1$ does not embed as a subgroup of any right-angled Artin group.
\end{cor}

Note that the proof of Theorem \ref{x2y2z2} shows that $\pi_1(S_{-1})$ is
not even {\it residually} isomorphic to a subgroup of a right-angled Artin group.
 

\section{The word and conjugacy problems for right-angled Artin groups}
\label{extraSect}


The word and conjugacy problems for right-angled Artin groups 
have been widely studied and are well understood
\cite{Green,HM,HW,Serv,VWyk,Wrath}. In this section we present
what we believe to be a particularly simple explanation of the solutions
as outlined in \cite{Serv}. We use, from the previous section, the notion of
a ``dissection of $D$ consistent with the boundary labelling'' which is equivalent
(in fact dual to) a Dehn/Van Kampen diagram. We note also that hypercurves of a
dissection associated to a generator $t$ play the same role as \emph{$t$-corridors}
used by various authors in the study of van Kampen diagrams.

Consider a right-angled Artin group $A_\Delta$
with defining graph $\Delta$. The word problem in $A_\Delta$ is solved by the
following ``algorithm'':
 
Let $w=u_1^{\ep_1}u_2^{\ep_2}\ldots u_k^{\ep_k}$, with $\ep_i\in\{\pm 1\}$, 
be a freely reduced word in the generators of $A_\Delta$. We say that $w$ is
\emph{$\Delta$-reduced} if there exist no $1\leq i < j\leq k$ such that
$u_i^{\ep_i}=(u_j^{\ep_j})^{-1}$
and $(u_i,u_r)$ is an edge of $\Delta$ for all $i< r<j$. Obviously,
if $w$ is not $\Delta$-reduced then it may easily be simplified by removing 
the letters $u_i^{\ep_i}$ and $u_j^{\ep_j}$
to give a shorter word representing the same element of $A_\Delta$.
This simplification process continues until one reaches a $\Delta$-reduced word.
We claim that the element of $A_\Delta$
under consideration is trivial if and only if this $\Delta$-reduced word is trivial.

The proof of this last claim is as follows:
suppose that $w$ is a nontrivial word that represents the identity
element in $A_\Delta$. Let $D$ be a disk whose boundary is subdivided into an
edge path labelled by the word $w$. Then, as in the proof of Claim A,
we can construct a dissection of $D$ which is consistent with the boundary
labelling. One can now choose a non-simple hypercurve $c$ 
(i.e: one which is not just a simple closed curve) which is
innermost in the following sense: any other non-simple
hypercurve which is disjoint from $c$ separates $c$ from
the vertex on $\partial D$ which marks where the word $w$ begins and ends.
It follows that the letters at the endpoints of $c$ commute with
all the letters between them in the word $w$ (since $c$ must cut the
corresponding hypercurves), and hence that $w$ is not $\Delta$-reduced. Thus
the only $\Delta$-reduced word which represents the identity is the trivial word.
 
The above result may be strengthened slightly, and extended to the case of 
cyclic words as in Proposition \ref{wordconj} below (see also \cite{Serv},
and the remark below).
A cyclically reduced word $w$ in the generators of $A_\Delta$
is said to be \emph{cyclically $\Delta$-reduced} if no cyclic permutation of $w$ 
contains a subword of the form $u^{\ep}v_1^{\ep_1}\ldots v_n^{\ep_n}u^{-\ep}$ where
$(u,v_i)$ is an edge of $\Delta$ for every $i=1,\ldots,n$. On the other hand,
two words in the generators of $A_\Delta$ are said to be 
\emph{$\Delta$-equivalent} if they are
related by a finite sequence of \emph{commutations} -- substitution of a 
subword $u^\ep v^\nu$ by $v^\nu u^\ep$ 
under the constraint that  $(u,v)$ be an edge of $\Delta$ -- and are said to be 
\emph{cyclically $\Delta$-equivalent} if they are related by a finite sequence 
of commutations
and cyclic permutations. Note that these equivalence operations do 
not change the length of a word.

\begin{prop}\label{wordconj}
Let $A_\Delta$ be the right-angled Artin group with defining graph $\Delta$.
\begin{itemize}
\item[\rm(i)] Two $\Delta$-reduced words represent the same element of $A_\Delta$ 
if and only if they are $\Delta$-equivalent. 
\item[\rm(ii)] Two cyclically $\Delta$-reduced words represent the same conjugacy 
class in $A_\Delta$ if and only if they are cyclically $\Delta$-equivalent.
\end{itemize}
\end{prop}   

\begin{proof}
(i)\qua Take $\Delta$-reduced words $x$ and $y$ representing the same element
in $A_\Delta$, and construct a dissection consistent with boundary
labelling for a disk $D$ with boundary labelled by the word $xy^{-1}$. 
Since $x$ and $y$ are $\Delta$-reduced we may suppose that every non-simple 
hypercurve joins a letter in $x$ with a similar letter in $y$. But then it is 
easy to see how to construct a finite sequence of commutations to show
that $x$ and $y$ are $\Delta$-equivalent, given that hypercurves which 
cross correspond to letters which commute.

(ii)\qua  Take $\Delta$-reduced words $x$ and $y$ representing 
conjugate elements of $A_\Delta$, so
$x=wyw^{-1}$ for some word $w$. Now construct a dissection consistent with boundary
labelling for a disk $D$ with boundary labelled by the word $x^{-1}wyw^{-1}$, 
and glue together
the two subsegments of $\partial D$ labelled by the word $w$ to give an annulus 
with boundaries
labelled by the words $x$ and $y$ respectively, and which carries a dissection 
which is consistent
with these boundary labels.  Again, since $x$ and $y$ are
$\Delta$-reduced, every non-simple hypercurve joins a letter in $x$ 
(on one side of the annulus) with
a letter in $y$ (on the other side), and one easily shows that $x$ and $y$
are cyclically $\Delta$-equivalent.
\end{proof}
\Remark The second part of Proposition \ref{wordconj} was 
given by Servatius in \cite{Serv} and yields the
polynomial time solution to the conjugacy problem alluded to in the Introduction.
The first part was stated without proof in \cite{Serv}. The result 
(or at least the analogous result stated in terms of ``expressions''
and ``syllables'' in Section \ref{Sect5} below)
was however treated in the more general context of graphs of groups by E. R. Green
in \cite{Green} (see also \cite{HM} for statements) and a different proof of
Green's result was given in \cite{HW}, using Dehn diagrams.
We claim that our proofs may be easily modified to apply also in the context of
graphs of groups. One simply allows each non-simple hypercurve of a
dissection to be, not just an arc, but rather a properly embedded trivalent graph
whose edges are transversely oriented and labelled by syllables
associated to a common vertex (non-trivial elements of a common
vertex group) in such a way that the product of labels read around
any vertex of the hypercurve is the identity. This branching allows for
``amalgamation'' and ``un-amalgamation'' of syllables in a given expression.
We leave the details to the interested reader.


\section{Right-angled Artin groups in pure surface braid groups}\label{Sect5}


The principal aim of this section is to prove the following:

\begin{thm}\label{thm9}
 Let $\Delta$ be a graph without loops or multiple edges, and let $n$
denote the number of vertices in $\Delta$. Then the right-angled Artin
group $A_\Delta$ associated to $\Delta$ can be faithfully represented 
as a subgroup of a pure surface braid group $PB_n(S)$,
where $S$ is a compact surface with nonempty boundary. 
 \end{thm}

Since pure surface braid groups of compact surfaces with boundary are
bi-orderable by \cite{G-M} we have as an immediate corollary that
right-angled Artin groups are bi-orderable. Bi-orderability of
right-angled Artin groups was previously proved by Duchamp and Thibon
\cite{DT} by generalising the Magnus construction for free groups.

Theorem \ref{thm9} will be eventually established by proving Proposition
\ref{prop11} below.  We fix some notation first.
Let $S$ denote 
a surface (not necessarily connected) which is obtained by plumbing together,
in an arbitrary manner, a collection $A_1,A_2,\ldots,A_n$ of disjoint annuli.
We make this more precise as follows.

For convenience we ``coordinatise" each annulus $A_i$ via a
homeomorphism to $\R/N\Z\times[0,1]$, where $N$ is a sufficiently
large natural number. By a \emph{patch} (on the annulus $A_i$) we
mean a closed subset of $A_i=\R/N\Z\times[0,1]$ of the form
$[t,t+1]\times[0,1]$ for some $t\in \R/N\Z$. The construction of
$S$ depends upon
 \begin{itemize}
 \item[\rm(1)] a choice of a finite collection $\Cal P$ of mutually disjoint
patches on the annuli $A_1,\ldots,A_n$;
 \item[\rm(2)] a choice of a family $\Cal B$ of ordered pairs $(Q,Q')$ where
$Q$ and $Q'$ are elements of $\Cal P$ coming from different annuli; and
 \item[\rm(3)] for each pair $(Q,Q')\in\Cal B$, a choice of homeomorphism
$h\co Q\to Q'$ such that, if $Q=[s,s+1]\times [0,1]\subset A_i$
and $Q'=[r,r+1]\times [0,1]\subset A_j$ say, then $h(x,y)=(f(y),g(x))$ 
where $f\co [0,1] \to [r,r+1]$ and $g: [s,s+1]\to [0,1]$ are homeomorphisms.
 \end{itemize}

\begin{defn}\label{defineS}
 The surface $S$ is the defined to be the
identification space
\[
S= 
\raise.75ex\hbox{$(\coprod\limits_{i=1}^n 
A_i)$}/\lower.75ex\hbox{$\approx$}\ ,
\]
 where $\approx$ is the equivalence relation determined by the family of
homeomorphisms $h\co Q\to Q'$ for $(Q,Q')\in\Cal B$. 

Given the surface $S$ defined in this way we define the graph $\Delta(S)$
as follows. The vertex set of $\Delta(S)$ is in one-to-one correspondance
with the annuli of $S$ and is written $\Cal V=\{a_1,a_2,\ldots,a_n\}$.
The pair $(a_i,a_j)$ represents an edge  of  $\Delta(S)$ if and only
if $A_i$ and $A_j$ are disjoint subsets of $S$.
\end{defn}

\begin{rem}\label{circlesremark} 
It is clear that every graph without loops or multiple edges may be realised as
$\Delta(S)$ for some surface constructed as above.
Note, however, that for a given graph the construction is far from unique.
In many cases a given graph may be realised by building a \emph{planar} surface $S$
(which results in an embedding of the associated right-angled Artin group in a
classical braid group -- see Section \ref{Sect6}). 
This construction requires drawing a 
family of simple closed curves in the plane which are mutually transverse
and have the desired intersection properties. In Figure \ref{F:circles} 
we illustrate such circle drawings which realise the $5$-cycle of
Figure \ref{F:0mod2}, the graph $\Delta(K_5)$, and the graph
$\Delta_1\cup\Delta_2$ of Figure \ref{F:3mod4} (and in particular its full
subgraph $\Delta_1=\Delta(K_{3,3})$).
The surface $S$ is obtained in each case by taking a planar neighbourhood
of the union of circles.  Note that the fact that the diagram of circles in Figure 
\ref{F:circles}(a) has the symmetry of an $5$-cycle is a consequence of the
fact that the $5$-cycle is self-dual. In general, one should expect a much more
complicated circle diagram when $\Delta$ is an $n$-cycle for $n>5$. 
\end{rem}

We now consider the $n$-string braid group on the
surface $S$ of Definition \ref{defineS}. It will be 
convenient to let $p_i$ denote the point $(0,0)$ and $q_i$ the point
$(0,\frac{1}{2})$ of the $i^{\rm th}$ annulus $A_i$, for each $i=1,\ldots,n$.
We may suppose in the above construction that these points
are not contained in any of the patches of $\Cal P$. We take as our set of
``marked points'' the set $\Cal Q=\{q_i:i=1,\ldots,n\}$. For $1\leq i \leq
n$, we define $\beta_i\in B_n(S)$ to be the braid that moves the marked
point $q_i$ once around $A_i$ in the positive direction, and leaves all
the other points in $\Cal Q$ fixed. More precisely, the braid $\beta_i$ is
represented by the path of $n$-tuples,
$\beta_i(t)=(q_1,\ldots,q_{i-1},\gamma_i(t),q_{i+1},\ldots,q_n)$, where
$\gamma_i(t)=(Nt,\frac{1}{2})\in A_i=\R/ N\Z\times [0,1]$ for $t\in[0,1]$.  
We observe that setting $\varphi(a_i)= \beta^2_i$, for each $i=1,\ldots,n$,
defines a homomorphism from $A_{\Delta(S)}$ to $PB_n(S)$. Theorem \ref{thm9} is
now proved by establishing the following:

\begin{figure}[ht!] 
\begin{center}
\begin{picture}(0,0)%
\includegraphics{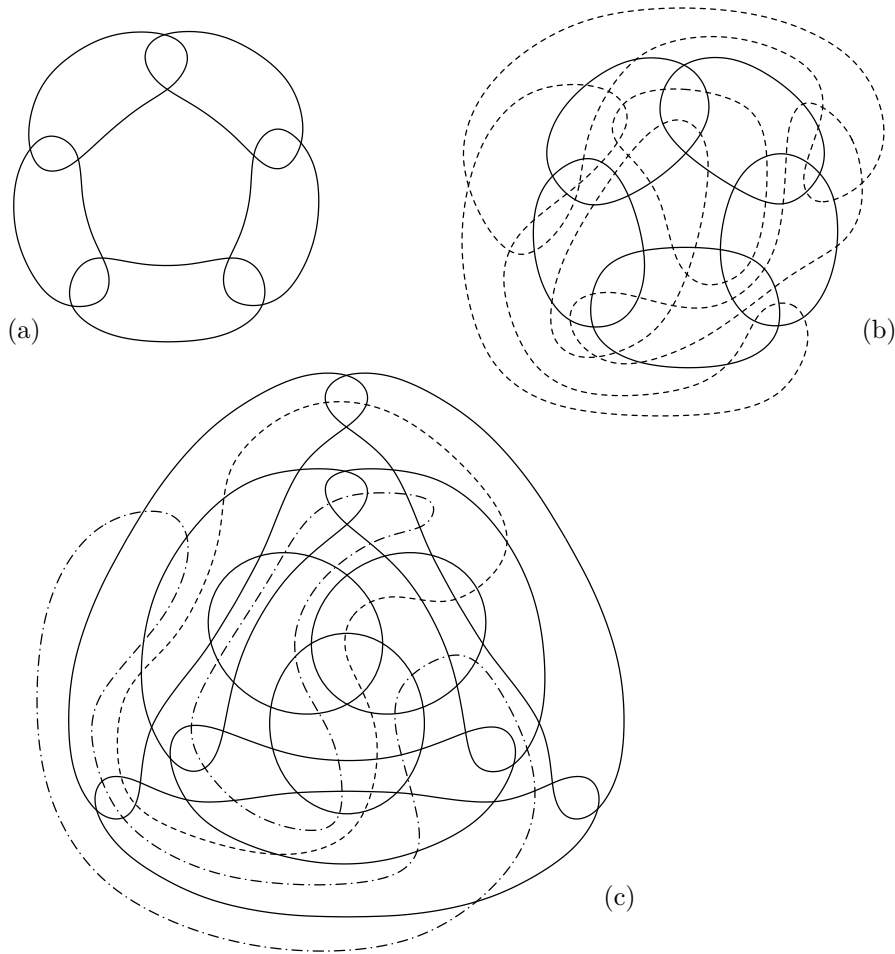}%
\end{picture}%
\setlength{\unitlength}{1066sp}%
\begingroup\makeatletter\ifx\SetFigFont\undefined%
\gdef\SetFigFont#1#2#3#4#5{%
  \reset@font\fontsize{#1}{#2pt}%
  \fontfamily{#3}\fontseries{#4}\fontshape{#5}%
  \selectfont}%
\fi\endgroup%
\begin{picture}(20396,22020)(4726,-23090)
\put(18601,-21961){\makebox(0,0)[lb]{\smash{\SetFigFont{10}{12.0}{\rmdefault}{\mddefault}{\updefault}(c)}}}
\put(4726,-8761){\makebox(0,0)[lb]{\smash{\SetFigFont{10}{12.0}{\rmdefault}{\mddefault}{\updefault}(a)}}}
\put(24601,-8761){\makebox(0,0)[lb]{\smash{\SetFigFont{10}{12.0}{\rmdefault}{\mddefault}{\updefault}(b)}}}
\end{picture}
\end{center}
\caption{Circle diagrams which yield planar surfaces $S$ with $\Delta(S)$ equal to
(a) the $5$-cycle; (b) $\Delta(K_5)$; and (c) $\Delta_1\cup\Delta_2$ with the
full subgraph $\Delta_1=\Delta(K_{3,3})$ being 
realised by the plain (non-dashed) circles.}\label{F:circles}
\end{figure}

\begin{prop}\label{prop11} Let $S$ denote any surface constructed from $n$ annuli
as in Definition \ref{defineS}.  Then the representation 
$\varphi\co A_{\Delta(S)} \to PB_n(S)$ defined by
$\varphi(a_i)= \beta^2_i$ is faithful.
\end{prop}

Let $F$ denote the $n$-punctured surface $S\setminus\Cal Q$. We note that
(pure) braids induce homeomorphisms of $F$ which are unique up to ambient
isotopy of $F$ and well-defined on braid equivalence classes. Thus we have
a homomorphism $\mu\co PB_n(S)\to \Cal M(F)$, where $\Cal M(F)$ denotes
the mapping class group of the surface $F$ relative to the boundary
$\partial F=\partial S$.

By a \emph{diagram in $F$} we mean a smooth embedding $E : \{
1,\ldots,n\}\times [0,1)\to F$ of $n$ disjoint half-open intervals into
$F$ which may be extended to a smooth embedding $E' \co \{
1,\ldots,n\}\times [0,1]\to S$ of $n$ disjoint closed intervals into $S$
by setting $E'(i,1)=q_i$. Let $E_i\co [0,1)\to F$ denote the $i$th curve
in the diagram $E$ in $F$, that is $E_i(t)=E(i,t)$. We shall say that two
diagrams in $F$ are \emph{equivalent} if there is an ambient isotopy of
$F$ relative to $\partial F$ which carries one diagram onto the other.
(Note that considering ambient isotopy classes relative to the boundary
means that the basepoints $p_i$ are always fixed while the open end of
each curve in a curve diagram is effectively anchored in the vicinity of
the associated puncture).

The group $PB_n(S)$ clearly acts on the set of equivalence classes of
diagrams in $F$ via the mapping class group (i.e: via the homomorphism
$\mu\co PB_n(S)\to \Cal M(F)$).

We define the diagram $E^0$ in $S$ to consist of the $n$ straight line
segments joining each point $p_i$ with the corresponding $q_i$. Thus
$E^0(i,t)=E^0_i(t)$ is the point $(0,\frac{t}{2})$ in $A_i$, for each
$i=1,\ldots,n$ and each $t\in[0,1]$.

\begin{figure}[ht!] 
\centerline{
\begin{picture}(0,0)%
\includegraphics{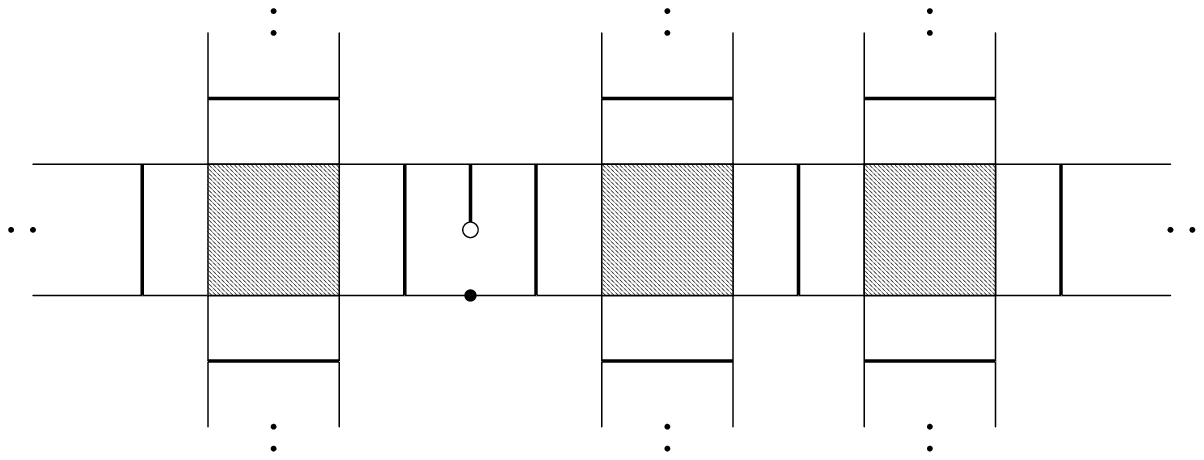}%
\end{picture}%
\setlength{\unitlength}{1381sp}%
\begingroup\makeatletter\ifx\SetFigFont\undefined%
\gdef\SetFigFont#1#2#3#4#5{%
  \reset@font\fontsize{#1}{#2pt}%
  \fontfamily{#3}\fontseries{#4}\fontshape{#5}%
  \selectfont}%
\fi\endgroup%
\begin{picture}(16330,6130)(836,-8526)
\put(1801,-5536){\makebox(0,0)[lb]{\smash{\SetFigFont{10}{12.0}{\rmdefault}{\mddefault}{\updefault}$A_i$}}}
\put(7351,-5386){\makebox(0,0)[lb]{\smash{\SetFigFont{10}{12.0}{\rmdefault}{\mddefault}{\updefault}$q_i$}}}
\put(7276,-6661){\makebox(0,0)[lb]{\smash{\SetFigFont{10}{12.0}{\rmdefault}{\mddefault}{\updefault}$p_i$}}}
\end{picture}
}
\caption{Cutting arcs on the annulus $A_i$}\label{F:cuttingarcs} 
\end{figure}

We now cut $F$ into simply connected pieces by introducing a family $\Cal
C$ of mutually disjoint intervals properly embedded
in $(F,\partial F)$ as indicated in Figure \ref{F:cuttingarcs}.
 Namely, we cut from $q_i$ to the point $(0,1)$ in $A_i$ (directly
opposite $p_i$), and we cut across the annulus $A_i$ once between every
pair of adjacent gluing patches, and once between the interval $\{
0\}\times [0,1]$ and each neighbouring gluing patch. In the case of
an annulus which is disjoint from all other annuli we make at least one
cut across the annulus in addition to the cut starting at $q_i$.
The elements of $\Cal C$ shall be called \emph{cutting arcs}.  A curve
diagram is said to be \emph{reduced} if it is in general position with
respect to the family $\Cal C$ of cutting arcs and its (geometric)
intersection number with any cutting arc is as small as possible within
its equivalence class.  It is not too hard to see that any diagram is
equivalent to a reduced one, and that this reduced representative is
unique up to an ambient isotopy through a family of reduced diagrams.
We may make this statement more precise as follows.

If we orient each cutting arc transversely, then each diagram $E$ in $F$
which is in general position with respect to $\Cal C$ may be associated
with an $n$-tuple $(w_1,\ldots,w_n)$ where each $w_i$ is a word in the
alphabet $\Cal C^\pm=\{ c,c^{-1} : c\in \Cal C\}$ obtained by reading off
with sign the sequence of crossings that the curve $E_i$ makes with $\Cal
C$ while travelling from the boundary point $p_i$ toward the marked point
$q_i$. We shall call $(w_1,\ldots,w_n)$ the \emph{crossing sequence} of the
diagram $E$, and refer to $w_i$ as the crossing sequence of $E_i$, for
each $i=1,\ldots,n$.
 A diagram $E$ is reduced if and only if each component of its crossing
sequence is a freely reduced word. Moreover, two reduced diagrams which
are equivalent have the same crossing sequence. This is simply because any
ambient isotopy of diagrams will induce a finite sequence of modifications
of the associated crossing sequences which consist only of trivial
insertions and deletions.

For each $i=1,\ldots,n$ we define the word $\alpha_i$ as follows. Let
$E^{(i)}$ denote a reduced diagram equivalent to $\beta_i(E^0)$. Then
$\alpha_i$ is defined as the crossing sequence for $E^{(i)}_i$.  In other
words, $\alpha_i$ is the crossing sequence for a reduced path which
travels exactly once around $A_i$ from $p_i$ to $q_i$.  Let $1\leq i\leq
n$. We say that a reduced diagram $E$ has an \emph{$i$-tail} if the
crossing sequence $w_i$ for $E_i$ ends with $\alpha_i$ or $\alpha_i^{-1}$.
More generally, we say that any diagram has an $i$-tail if it is
equivalent to a reduced diagram with an $i$-tail. By our previous remarks,
this notion is well defined on equivalence classes of diagrams in $F$.

\begin{lemma}\label{lemm1} Let $E$ be a diagram in $F$, 
and let $k\in\Z\setminus\{ 0\}$.
 \begin{itemize}
 \item[\rm(i)] If $E$ has both an $i$-tail and a $j$-tail, for $i\neq j$, 
then $(a_i,a_j)\in\Cal E$, i.e.~the generators $a_i$ and $a_j$ of 
$A_{\Delta(S)}$ commute.
 \item[\rm(ii)] If $E$ does not have an $i$-tail, then $\beta_i^{2k}(E)$ does 
(have an $i$-tail).
 \item[\rm(iii)] Let $(a_i,a_j)\in\Cal E$. Then $E$ has a $j$-tail if and 
only if  $\beta_i^{2k}(E)$ has a $j$-tail.
 \end{itemize} 
 \end{lemma} 

\begin{proof} Without loss of generality, we may suppose that all diagrams 
considered are reduced.

(i)\qua The (reduced) diagram $E$ having an $i$-tail means that the curve
$E_i$ makes a complete circuit along the length of the annulus $A_i$ thus
crossing transversely every annulus $A_k$ which intersects $A_i$ (in fact,
crossing every gluing patch). Since the different curves $E_i$ and $E_j$
of the diagram $E$ are disjoint, it is clear that $E$ cannot therefore
have both an $i$-tail and a $j$-tail unless $A_i$ and $A_j$ are disjoint,
in which case $(a_i,a_j)\in\Cal E$.

We leave (ii) and (iii) as exercises for the reader. Note that in part (ii) 
at least a double twist $\beta^2_i$ is needed to ensure the appearance of
an $i$-tail.
\end{proof}

To prove Proposition \ref{prop11}, we consider a nontrivial element $x\in
A_\Delta$. Each nontrivial element can by represented by a nontrivial
``reduced expression'' as defined below.

An \emph{expression} for $x\in A_\Delta$ of \emph{syllable length} $\ell$
is a word in the generators $\{a_1,\ldots,a_n\}$ representing $x$ and of the
form
 \[
 X=a_{i_1}^{k_1}\cdots a_{i_\ell}^{k_\ell}\,\qquad \text{ where }
k_j\in\Z\setminus\{ 0\}\text{ and } i_j\in\{1,\ldots,n\} \text{ for } 1\leq
j\leq \ell\,.
 \]
 Briefly, an amalgamation is simply lumping together adjacent syllables of
the same type to shorten the syllable length of $X$, and shuffle
equivalence of expressions is generated by exchanges of adjacent syllables
$a_i^k$, $a_j^m$ where $(a_i,a_j)\in\Cal E$. These operations yield new
expressions for the same element of $A_\Delta$.  A \emph{reduced
expression} is an expression which is not shuffle equivalent to one which
admits an amalgamation. The key properties of reduced expressions are that
(i) every element $x\in A_\Delta$ may be represented by a reduced
expression $X$ (obtained from a given expression by a finite process of
shuffling and amalgamating), and (ii) any two reduced expressions $X,X'$
for the same element $x\in A_\Delta$ are shuffle equivalent. We refer to
\cite{Green} or \cite{HM} for more details -- see also \cite{CP}.
Note that property (ii) is a consequence of the following Lemma.

We say that a (reduced) expression \emph{ends in $a_i$} if it is shuffle
equivalent to an expression of the form $a_{i_1}^{k_1}\cdots
a_{i_\ell}^{k_\ell}$ where $i_1=i$.

\begin{lemma}\label{lemm12}
 Let $X$ be a reduced expression. Then $X$ ends in $a_i$ if and only if
the diagram $\varphi(X)(E^0)$ has an $i$-tail. 
\end{lemma}

\begin{proof} We proceed by induction of $\ell(X)$.  Note firstly that, for all
$i=1,\ldots,n$, the trivial expression $X=1$ does not end in $a_i$ while at
the same time $E^0=\varphi(1)(E^0)$ does not have an $i$-tail.

Suppose $X$ ends in $a_i$, then $X=a_i^kX'$ for some $k\in\Z\setminus\{
0\}$ and some reduced expression $X'$ where $\ell(X')=\ell(X)-1$ and $X'$
does not end in $a_i$.  By induction $E=\varphi(X')(E^0)$ has no
$i$-tail. Thus, by Lemma \ref{lemm1}(ii), $\varphi(X)(E^0)\sim
\beta^{2k}_i(E)$ does have an $i$-tail, as required.
 
Now, suppose that $\varphi(X)(E^0)$ has a $j$-tail. We shall show that $X$
ends in $a_j$. Since $X$ is necessarily nontrivial there exists an $i$
for which $X$ ends in $a_i$, and we may assume that $i\neq j$. By the
argument just given we may deduce that $\varphi(X)(E^0)$ has an $i$-tail.
But now by Lemma \ref{lemm1}(i) we have $(a_i,a_j)\in\Cal E$. Writing
$X=a_i^kX'$ for some $k\in\Z\setminus\{ 0\}$ and some reduced expression
$X'$ with $\ell(X')=\ell(X)-1$, we have
$\varphi(X)(E^0)=\beta_i^{2k}(\varphi(X')(E^0))$ and by Lemma
\ref{lemm1}(iii) we deduce that $\varphi(X')(E^0)$ also has a $j$-tail.
Now, by induction, $X'$ ends in $a_j$ and, since the letters $a_i$ and
$a_j$ commute, a single shuffle shows that $X$ also ends in $a_j$.  
\end{proof}

\begin{proof}[Proof of Proposition \ref{prop11}]
 Recall that the pure braid group $PB_n(S)$ acts on the punctured surface
by isotopy classes of homeomorphisms via the homomorphism $\mu\co
PB_n(S)\to \Cal M(F)$. Let $\varphi'$ denote the composition
$\mu\circ\varphi$. Our aim is now to prove that $\varphi'$ is injective. 
Note that any nontrivial element $x\in A_\Delta$ is represented by a nontrivial
reduced expression $X$, which ends in $a_i$ say. It follows from Lemma
\ref{lemm12} that $\varphi'(x)(E^0)$ has an $i$-tail, and hence is not
equivalent to the diagram $E^0$. That is to say that $\varphi'(x)$ is a
nontrivial mapping class. Therefore $\varphi'$ is injective. 
\end{proof}

\begin{rem}\label{realize}
 Note that, in order to show that the group $A_\Delta$ embeds in the pure
$n$-string braid group of $S$ we have constructed a faithful
representation $\varphi'\co A_\Delta\to\Cal M(F)$ of $A_\Delta$ in the
mapping class group of the $n$-punctured surface $S$ such that $\varphi'$
factors through $PB_n(S)$.  In fact it is not too hard to see that the
image of $\varphi'$ is actually realizable in the sense of Nielsen and
Kerckhoff \cite{Kerck}:  ${\rm im}(\varphi')$ can be realised by a group
of diffeomorphisms of $S$ which fix $\partial S$ and the distinguished
points $q_1,\ldots, q_n$. In other words $\varphi'$ factors through an
(injective) homomorphism $f\co A_\Delta\to\text{Diffeo}(S;\partial S,
q_1\ldots q_n)$. This is because, for each $i=1,2,\ldots,n$, the mapping class
$\varphi'(a_i)=\mu(\beta_i^2)$ may be represented by a diffeomorphism
$\alpha_i$ which is the identity outside the annulus $A_i$. In the case
that $(a_i,a_j)\in\Cal E$, the annuli $A_i$ and $A_j$ are disjoint subsets
of $S$ and it is clear that $\alpha_i$ and $\alpha_j$ commute as
diffeomorphisms. Thus all of the relations in the standard presentation of
$A_\Delta$ are satisfied in the group $\text{Diffeo}(S;\partial S,
q_1\ldots q_n)$. \end{rem}


\section{Embeddings in the classical braid groups}\label{Sect6}


So far we have developed a method for embedding right-angled Artin groups
in surface braid groups. It seems a much harder task to find embeddings in
classical braid groups. It is nevertheless conceivable that this might be
possible for all right-angled Artin groups, and in this
section we present some partial results in this direction.

Given a finite simple graph $\Delta$, we define the \emph{dual} 
or \emph{opposite} graph
$\Delta^{\rm op}$ to be the finite simple graph with the same 
vertex set as $\Delta$ where
$(u,v)$ is defined to be an edge of $\Delta^{\rm op}$ precisely when
$u$ and $v$ are distinct non-adjacent vertices of $\Delta$.

\begin{prop}\label{planar} 
Let $A_\Delta$ be the right-angled Artin group defined by the graph
$\Delta$. If the opposite graph $\Delta^{\rm op}$ is a planar graph then
the $A_\Delta$ embeds as a subgroup of the classical pure braid group
$PB_l$ for sufficiently large $l\in \N$, and also in ${\rm
Diff}(D^2;\partial D^2, p_1\ldots p_l)$. 
\end{prop}

\begin{figure}[ht!]
\begin{center}
\includegraphics[width=13cm]{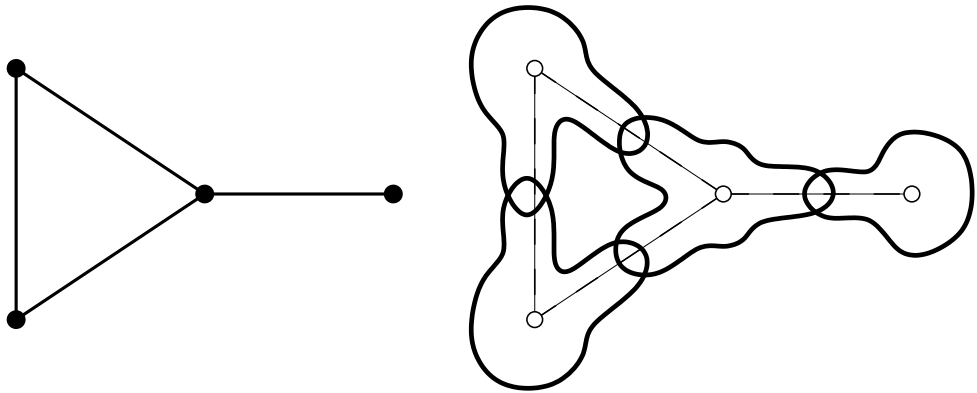}
\end{center}
\caption{How to construct a circle diagram for a planar
surface $S$ such that $\Delta(S)=\Delta$, starting from
the planar dual graph $\Delta^{{\rm op}}$ (shown on the left).}\label{HtoS}
\end{figure}

\begin{proof} One can use the construction method suggested in Figure \ref{HtoS}
to obtain a planar surface $S$ such that $\Delta(S)=\Delta$.
By Theorem \ref{thm9} we obtain an embedding of $A_\Delta$ in
the pure braid group, and the diffeomorphism
group, of $S$. By glueing a twice-punctured disk to all but one of the
components of $\partial S$ and extending all automorphisms of $S$ by the
identity on these twice-punctured disks, we get faithful representations
in the the braid groups and the diffeomorphism groups, respectively, of an
$l$-punctured disk (see Theorem 2.3 of \cite{PR1}, also \cite{PR2}). \end{proof}

\begin{figure}[ht!] %
\begin{center}
\begin{picture}(0,0)%
\includegraphics{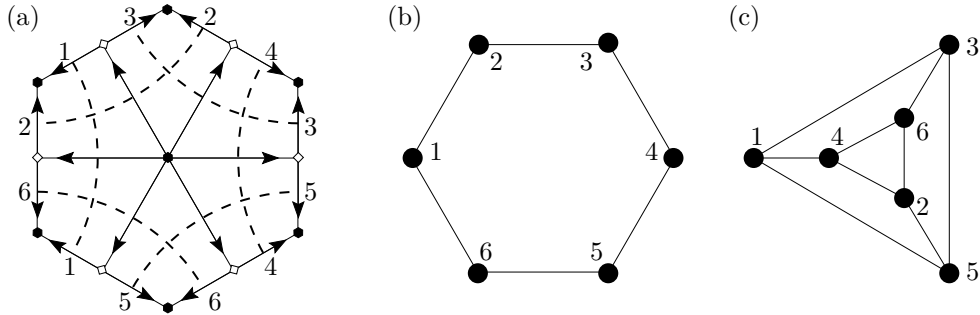}%
\end{picture}%
\setlength{\unitlength}{1658sp}%
\begingroup\makeatletter\ifx\SetFigFont\undefined%
\gdef\SetFigFont#1#2#3#4#5{%
  \reset@font\fontsize{#1}{#2pt}%
  \fontfamily{#3}\fontseries{#4}\fontshape{#5}%
  \selectfont}%
\fi\endgroup%
\begin{picture}(14637,4685)(601,-3891)
\put(4442,-69){\makebox(0,0)[lb]{\smash{\SetFigFont{10}{12.0}{\rmdefault}{\mddefault}{\updefault}$4$}}}
\put(5026,-1186){\makebox(0,0)[lb]{\smash{\SetFigFont{10}{12.0}{\rmdefault}{\mddefault}{\updefault}$3$}}}
\put(5026,-2161){\makebox(0,0)[lb]{\smash{\SetFigFont{10}{12.0}{\rmdefault}{\mddefault}{\updefault}$5$}}}
\put(4426,-3286){\makebox(0,0)[lb]{\smash{\SetFigFont{10}{12.0}{\rmdefault}{\mddefault}{\updefault}$4$}}}
\put(3601,-3811){\makebox(0,0)[lb]{\smash{\SetFigFont{10}{12.0}{\rmdefault}{\mddefault}{\updefault}$6$}}}
\put(2251,-3811){\makebox(0,0)[lb]{\smash{\SetFigFont{10}{12.0}{\rmdefault}{\mddefault}{\updefault}$5$}}}
\put(1426,-3286){\makebox(0,0)[lb]{\smash{\SetFigFont{10}{12.0}{\rmdefault}{\mddefault}{\updefault}$1$}}}
\put(751,-2161){\makebox(0,0)[lb]{\smash{\SetFigFont{10}{12.0}{\rmdefault}{\mddefault}{\updefault}$6$}}}
\put(751,-1186){\makebox(0,0)[lb]{\smash{\SetFigFont{10}{12.0}{\rmdefault}{\mddefault}{\updefault}$2$}}}
\put(1351,-61){\makebox(0,0)[lb]{\smash{\SetFigFont{10}{12.0}{\rmdefault}{\mddefault}{\updefault}$1$}}}
\put(2326,464){\makebox(0,0)[lb]{\smash{\SetFigFont{10}{12.0}{\rmdefault}{\mddefault}{\updefault}$3$}}}
\put(3526,464){\makebox(0,0)[lb]{\smash{\SetFigFont{10}{12.0}{\rmdefault}{\mddefault}{\updefault}$2$}}}
\put(601,464){\makebox(0,0)[lb]{\smash{\SetFigFont{10}{12.0}{\rmdefault}{\mddefault}{\updefault}(a)}}}
\put(6901,-1561){\makebox(0,0)[lb]{\smash{\SetFigFont{10}{12.0}{\rmdefault}{\mddefault}{\updefault}$1$}}}
\put(10126,-1561){\makebox(0,0)[lb]{\smash{\SetFigFont{10}{12.0}{\rmdefault}{\mddefault}{\updefault}$4$}}}
\put(7801,-211){\makebox(0,0)[lb]{\smash{\SetFigFont{10}{12.0}{\rmdefault}{\mddefault}{\updefault}$2$}}}
\put(9151,-211){\makebox(0,0)[lb]{\smash{\SetFigFont{10}{12.0}{\rmdefault}{\mddefault}{\updefault}$3$}}}
\put(7651,-3061){\makebox(0,0)[lb]{\smash{\SetFigFont{10}{12.0}{\rmdefault}{\mddefault}{\updefault}$6$}}}
\put(9301,-3061){\makebox(0,0)[lb]{\smash{\SetFigFont{10}{12.0}{\rmdefault}{\mddefault}{\updefault}$5$}}}
\put(6301,464){\makebox(0,0)[lb]{\smash{\SetFigFont{10}{12.0}{\rmdefault}{\mddefault}{\updefault}(b)}}}
\put(11401,464){\makebox(0,0)[lb]{\smash{\SetFigFont{10}{12.0}{\rmdefault}{\mddefault}{\updefault}(c)}}}
\put(11701,-1336){\makebox(0,0)[lb]{\smash{\SetFigFont{10}{12.0}{\rmdefault}{\mddefault}{\updefault}$1$}}}
\put(14926, 14){\makebox(0,0)[lb]{\smash{\SetFigFont{10}{12.0}{\rmdefault}{\mddefault}{\updefault}$3$}}}
\put(14926,-3361){\makebox(0,0)[lb]{\smash{\SetFigFont{10}{12.0}{\rmdefault}{\mddefault}{\updefault}$5$}}}
\put(14176,-2386){\makebox(0,0)[lb]{\smash{\SetFigFont{10}{12.0}{\rmdefault}{\mddefault}{\updefault}$2$}}}
\put(14176,-1261){\makebox(0,0)[lb]{\smash{\SetFigFont{10}{12.0}{\rmdefault}{\mddefault}{\updefault}$6$}}}
\put(12901,-1336){\makebox(0,0)[lb]{\smash{\SetFigFont{10}{12.0}{\rmdefault}{\mddefault}{\updefault}$4$}}}
\end{picture}
\end{center}
\caption{(a) A cubing of the orientable surface of genus 2, with the six
midplanes marked in dashed lines. The surface is obtained by identifying
edges with the same label. (b) The graph $\Delta$ symbolizing the
right-angled Artin group in which $\pi_1(S)$ can be embedded. (c) The
opposite graph $\Delta^{{\rm op}}$ is planar.}
\label{planarsurf}
\end{figure}

\begin{exam}\rm Let $\Delta$ be the hexagonal graph shown in
Figure \ref{planarsurf}(b).
The right-angled Artin group $A_\Delta$ contains the
fundamental group of the orientable surface of genus 2 (see Figure
\ref{planarsurf}(a)) which is embedded by the technique of Section \ref{Sect4}.
The opposite graph $\Delta^{{\rm op}}$ (shown in Figure \ref{planarsurf}(c))
is planar and therefore, by Proposition \ref{planar},
the group $A_\Delta$ (together with the surface subgroup)
embeds in a classical pure braid group $PB_l$ ($l=44$ suffices).
\end{exam}

\begin{prop}\label{planarcov}
If the graph $\Delta^{{\rm op}}$ has a
finite-sheeted covering which is planar, then the right-angled Artin group
symbolised by $\Delta$ embeds in the classical pure braid group $PB_l$ for
sufficiently large $l$, and also in 
$\hbox{Diff}(D^2;\partial D^2, p_1\ldots p_l)$. 
\end{prop}

\begin{proof}
 Let $\widetilde{\Delta}^{{\rm op}}$ be the planar cover of $\Delta^{{\rm
op}}$ (with say $k$ sheets). We define $\widetilde{A}_\Delta$ to be the
right-angled Artin group whose defining graph is the opposite graph of
$\widetilde{\Delta}^{{\rm op}}$. There is an obvious homomorphism $j \co
A_\Delta \to \widetilde{A}_\Delta$, given by sending a generator $h$ of
$A_\Delta$ to the $k$-letter word $\widetilde{h}^{(1)}\ldots
\widetilde{h}^{(k)}$, where
$\widetilde{h}^{(1)},\ldots,\widetilde{h}^{(k)}$ are the generators of
$\widetilde{A}_\Delta$ covering $h$. Note that
$\widetilde{h}^{(1)},\ldots,\widetilde{h}^{(k)}$ commute in
$\widetilde{A}_\Delta$.

We claim that the homomorphism $j \co A_\Delta \to \widetilde{A}_\Delta$
is injective. To see this, we consider a word $w$ in the letters $h_1^{\pm
1},\ldots, h_m^{\pm 1}$ which is a representative of the smallest possible
length of $[w]\in A_\Delta$. Equivalently, by Proposition \ref{wordconj},
 $w$ is $\Delta$-reduced: it has no subword of the form
$h_i w' h_i^{-1}$ (or $h_i^{-1} w' h_i$), where $i\in \{1,\ldots,m\}$, and
$w'$ is a word consisting entirely of letters that commute with $h_i$. But
then the ``obvious'' representative of $j([w])$ (given by the
letter-by-letter substitution $h_i\mapsto
\widetilde{h}^{(1)}_i\ldots\widetilde{h}^{(k)}_i$)  satisfies the same
condition. It follows that this word is also of the smallest possible
length, and in particular $j([w])\neq 0$. (See \cite{CP} for a similar
argument in a slightly different situation).

Thus $A_\Delta$ embeds in $\widetilde{A}_\Delta$, which in turn embeds in
$B_l$ and $\hbox{Diff}(D^2;\partial D^2, p_1\ldots p_l)$ for sufficiently
large $l$, by Proposition \ref{planar}.
\end{proof}

\begin{exam}\rm 
 The graph $\Delta(K_{3,3})$ shown in Figure \ref{F:chi-3} is self-dual. 
That is, $\Delta(K_{3,3})^{\rm op}\cong \Delta(K_{3,3})$.
This graph is not planar, but it has a planar double cover
(since it embeds in a M\"obius band).
This implies that the fundamental group of the nonorientable
surface of Euler characteristic $-3$ embeds in a pure classical braid group,
using Proposition \ref{planarcov}.  
\end{exam}

There are however many examples of right-angled Artin groups (notably the
groups $\Z^n$) which embed in $PB_l$, for some $l$,
but for which the above propositions
either do not apply, or do not provide the best strategy.
It might be that in many individual cases the direct approach of drawing 
circles in the plane, such as is illustrated
in Figure \ref{F:circles} (see Remark \ref{circlesremark}), is the most efficient.
These drawings yield embeddings in planar surface braid groups and as explained
in the proof of Proposition \ref{planar}
the planar surfaces may be embedded in multi-punctured disks in order to obtain
embeddings in $PB_l$ for a suitable $l$.

\begin{figure}[ht!] 
\centerline{
\begin{picture}(0,0)%
\includegraphics{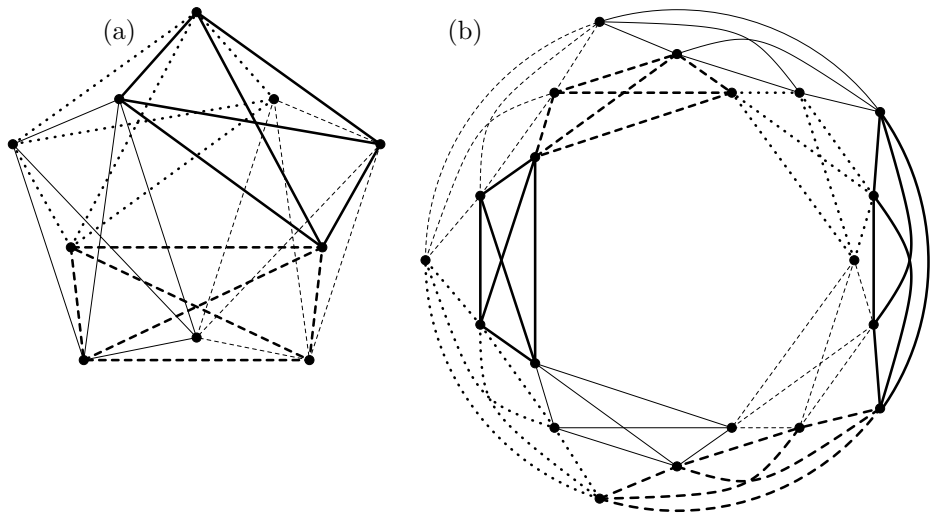}%
\end{picture}%
\setlength{\unitlength}{1066sp}%
\begingroup\makeatletter\ifx\SetFigFont\undefined%
\gdef\SetFigFont#1#2#3#4#5{%
  \reset@font\fontsize{#1}{#2pt}%
  \fontfamily{#3}\fontseries{#4}\fontshape{#5}%
  \selectfont}%
\fi\endgroup%
\begin{picture}(21453,11759)(637,-11407)
\put(2851,-361){\makebox(0,0)[lb]{\smash{\SetFigFont{10}{12.0}{\rmdefault}{\mddefault}{\updefault}(a)}}}
\put(10876,-361){\makebox(0,0)[lb]{\smash{\SetFigFont{10}{12.0}{\rmdefault}{\mddefault}{\updefault}(b)}}}
\end{picture}
}
\caption{(a) The graph $\Delta(K_5)^{{\rm op}}$; and (b) a double cover
$\widetilde{\Delta}^{{\rm op}}$ of the graph in (a).}
\label{Api1N}
\end{figure}

\begin{exam}\rm 
Consider $\Delta=\Delta(K_5)$. Recall that $A_{\Delta(K_5)}$ contains
an embedded copy of the graph braid group $RB_2(K_5)$ which is
isomorphic to the fundamental group of the closed surface of Euler
characteristic $-5$.  The circle diagram of Figure \ref{F:circles}(c)
yields an embedding of $A_\Delta$ in a planar surface braid group, and
hence in a classical pure braid group. This example may alternatively
be treated by a combination of ideas, as follows.  The graph
$\Delta^{\rm op}$ is shown in Figure \ref{Api1N}(a). Passing to the
double cover $\widetilde{\Delta}^{{\rm op}}$ shown in Figure
\ref{Api1N}(b), it suffices by the argument of Proposition
\ref{planarcov}, to consider the right-angled Artin group whose
defining graph is dual to $\widetilde{\Delta}^{{\rm op}}$. It is now a
relatively easy matter to produce an appropriate circle diagram by
modifying the idea of Figure \ref{HtoS} as suggested by Figure
\ref{closeup}, in order to embed the latter group in a classical pure
braid group.
 
Note that this, and every embedding discussed in this section,
factors through a group of actual diffeomorphisms of the $l$ times
punctured disk, as explained in Remark \ref{realize} of the previous section.

\begin{figure}[ht!] %
\begin{center}
\includegraphics[width=12cm]{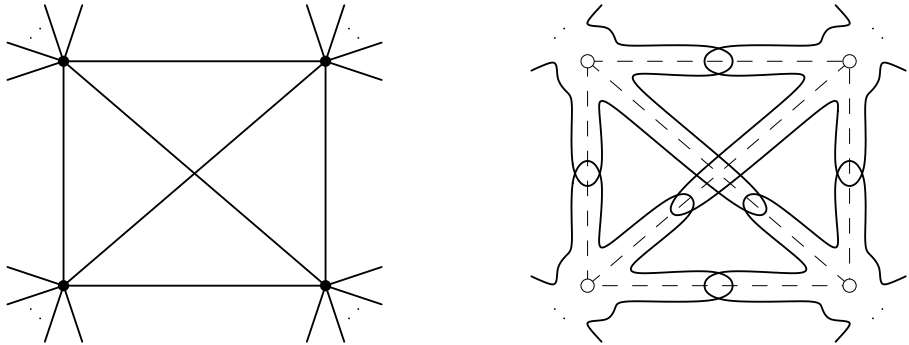}
\end{center}
\caption{How to modify the construction of the circle diagram 
for the surface $S$ at each of the
self-intersection points in the planar projection of the graph 
$\widetilde{\Delta^{\rm op}}$.}
\label{closeup}
\end{figure}
\end{exam}

The circle diagrams
in Figure \ref{F:circles} allow one to realise embeddings in a certain
$PB_l$ for each of the right-angled Artin groups used
in Section \ref{Sect4}. As a consequence, we have
the following:

\begin{prop}
With the exception of $RP^2$, $Kl$ and $S_{-1}$, the fundamental group of
every surface embeds in a pure braid group $PB_l$ for some $l$. 
\end{prop}

Observe that the group $\pi_1(RP^2)$ of order $2$ never embeds (since the braid
groups are torsion free), and $\pi_1(Kl)$ cannot embed in a pure braid group 
since the latter are bi-orderable while $\pi_1(Kl)$ is not
(it has generalised torsion). However, we do not know whether $\pi_1(Kl)$
embeds in $B_l$ for some $l$, or whether
$\pi_1(S_{-1})$ embeds in either $B_l$ or $PB_l$ for some $l$.

Finally, we make the following observation:

\begin{prop}\label{RBinbraid}
Let $G$ denote a finite graph. If $G$ or some finite index cover of $G$
is a planar graph, then, for any $n\geq 2$, the reduced braid group
$RB_n(G)$  embeds in a pure braid group $PB_l$ for some $l$.
\end{prop}

\begin{proof}
If $G$ itself is planar, then $\Delta(G)^{\rm op}$ is also planar and may
be constructed (in the plane) as follows: place a vertex of $\Delta(G)^{\rm op}$
at the midpoint of each edge of $G$, and an edge between two such vertices
if the corresponding edges of $G$ share a common vertex in $G$.
Note also that, if $G\to G_0$ is a finite index covering map (with $G$ planar),
then there is a naturally induced covering map
$\Delta(G)^{\rm op}\to \Delta(G_0)^{\rm op}$.
The result now follows by Proposition \ref{planarcov}.
\end{proof}

We do not know whether every reduced braid group may be embedded 
in a pure braid group $PB_l$. 
In connection with Propositions \ref{planarcov} and \ref{RBinbraid}
we note that it seems to be an open question
(or conjecture) as to whether every graph which has a planar
finite cover has in fact a planar double cover (or is
already planar), attributed to Henry Glover -- see
Bestvina's ``Questions in Geometric Group Theory''
problem list.   

Note, finally, that since $K_{3,3}$ and $K_5$ both embed in a M\"obius band
and therefore have planar double covers, Proposition \ref{RBinbraid}
gives yet another
method for embedding the closed Euler characteristic $-3$ and $-5$ surface
groups discussed above.


\Addresses\recd


\begin{thebibliography}{999}


\itemsep 1pt plus 1pt

\bibitem{Ab1}
{\bf A Abrams,} 
{\it Configuration Spaces and Braid Groups of Graphs},
PhD thesis, University of California at Berkeley (2000)

\bibitem{Ab2}
{\bf A Abrams,}
{\it Configuration spaces of colored graphs},
Geometriae Dedicata, to appear

\bibitem{AbGh}
{\bf A Abrams, R Ghrist,}
{\it Finding topology in a factory: configuration spaces},
Amer. Math. Monthly  109 no. 2 (2002) 140--150 

\bibitem{BB}
{\bf M Bestvina, N Brady,}
{\it Morse theory and finiteness properties of groups},  
Invent. Math. 129 (1997), no. 3, 445--470 

\bibitem{BH}
{\bf M\,R Bridson, A Haefliger,}
{\it Metric Spaces of Non-Positive Curvature},
Springer-Verlag (1999)

\bibitem{Ch}
{\bf R Charney,}
{\it The Tits conjecture for locally reducible Artin groups},
 Internat. J. Algebra Comput.  10 no. 6 (2000)  783--797

\bibitem{CP}
{\bf J Crisp, L Paris,} 
{\it The solution to a conjecture of Tits
on the subgroup generated by the squares of the generators of an Artin group},
Invent.~Math. 145 (2001) 19-36

\bibitem{DSS}
{\bf C Droms, B Servatius, H Servatius,}
{\it Surface subgroups of graph groups}, 
 Proc. Amer. Math. Soc.  106  (1989),  no. 3, 573--578

\bibitem{DT}
{\bf G Duchamp, J-Y Thibon,}
{\it Simple orderings for free partially commutative groups},
Internat. J. Algebra Comput. 2 (1992), no. 3, 351--355


\bibitem{G-M}
{\bf J Gonzalez-Meneses,}
{\it Ordering pure braid groups on compact, connected surfaces}, 
Pacific J.\ Math. 203 (2002), 369--378

\bibitem{GLR}
{\bf C Gordon, D\,D Long, A\,W Reid,} 
{\it Surface subgroups of Coxeter and Artin groups}. 
J. Pure Appl. Algebra 189 (2004), no. 1-3, 135--148.

\bibitem{Green}
{\bf E\,R Green,}
{\it Graph Products of Groups,}
Thesis, The University of Leeds, 1990.

\bibitem{Gro}
{\bf M Gromov,}
{\it Hyperbolic Groups}, Essays in Group Theory (S M Gersten, ed),
Springer-Verlag, MSRI Publ. 8 (1987) 75--263

\bibitem{HM}
{\bf S Hermiller, J Meier,}
{\it Algorithms and geometry for graph products of groups,}
J. Algebra 171 (1995), no. 1, 230--257

\bibitem{HW}
{\bf T Hsu, D\,T Wise,} 
{\it On linear and residual properties of graph products,}
Michigan Math.~J., to appear.

\bibitem{Hum}
{\bf S\,P Humphries,}
{\it On representations of Artin groups and the Tits conjecture,} 
J.~Algebra  169  (1994),  no.~3, 847--862

\bibitem{Kerck}
{\bf S Kerckhoff,}
{\it The Nielsen realisation problem},
Ann.~of Math.~(2) 117 (1983) 235--265


\bibitem{Lyn}
{\bf R\,C Lyndon,}
{\it The equation $a^2b^2=c^2$ in free groups},
Michigan Math.~J.~6 (1959) 155--164

\bibitem{LSbook}
{\bf R\,C Lyndon, P\,E Schupp,}
{\it Combinatorial Group Theory,}
Springer-Verlag, Berlin (1977).

\bibitem{LS}
{\bf R\,C Lyndon, M\,P Sch\"utzenberger,}
{\it The equation $a^M=b^Nc^P$ in a free group},
Michigan Math.~J.~9 (1962) 289--298

\bibitem{NR} 
{\bf G Niblo, L Reeves,}
{\it The geometry of cube complexes and the complexity of their fundamental 
groups},
Topology 37 (1998) 621--633

\bibitem{PR1}
{\bf L Paris, D Rolfsen,}
{\it Geometric subgroups of surface braid groups,}
Ann. Inst. Fourier 49 (1999),101-156. 

 \bibitem{PR2}
{\bf L Paris, D Rolfsen,}
{\it Geometric subgroups of mapping class groups,}
J.~Reine  Angew.~Math.~521 (2000),47-83.
 

\bibitem{RolfW}
{\bf D Rolfsen, B Wiest,}
{\it Free group automorphisms, invariant orderings and topological 
applications},
\agtref1{2001}{15}{311}{320}

\bibitem{Serv}
{\bf H Servatius,}
{\it Automorphisms of graph groups}, 
J.~Algebra 126 (1989) 34--60


\bibitem{VWyk}
{\bf L VanWyk,} 
{\it Graph groups are biautomatic},
J. Pure Appl. Algebra 94 (1994), no. 3, 341--352

\bibitem{Wrath}
{\bf C Wrathall,}
{\it Free partially commutative groups},
Combinatorics, computing and complexity (Tianjing and Beijing, 1988), 
195--216, Math.~Appl.~(Chinese Ser.), 1,
Kluwer Acad.~Publ., Dordrecht, 1989

\end{thebibliography}
\end{document}